\documentclass[letterpaper,12pt]{article}

\pdfoutput=1
\usepackage{graphicx}
\usepackage{amsmath}
\usepackage{amsfonts}
\usepackage{hyperref}
\usepackage{bm}

\makeatletter
\renewcommand\section{\@startsection {section}{1}{\z@}%
{-3.5ex \@plus -1ex \@minus -0.2ex}%
{2.3ex \@plus 0.2ex}%
{\normalfont\normalsize\bfseries}}

\renewcommand\subsection{\@startsection{subsection}{2}{\z@}%
{-3.25ex \@plus -1ex \@minus -0.2ex}%
{1.5ex \@plus 0.2ex}%
{\normalfont\normalsize\bfseries}}

\def\@seccntformat#1{\csname the#1\endcsname.\quad}
\makeatother

\newcommand\psiprime{\psi\hspace{0.05em}^{\prime}}

\begin{document}

\setlength{\baselineskip}{4.5ex}

\noindent
{\LARGE \bf Multivariate subjective fiducial inference}\\[3ex]

\noindent
{\bf Russell J. Bowater}\\
\emph{Independent researcher, Sartre 47, Acatlima, Huajuapan de Le\'{o}n, Oaxaca, C.P.\ 69004,
Mexico. Email address: as given on arXiv.org. Twitter profile:
\href{https://twitter.com/naked_statist}{@naked\_statist}\\ Personal website:
\href{https://sites.google.com/site/bowaterfospage}{sites.google.com/site/bowaterfospage}}
\\[2ex]

\noindent
{\small \bf Abstract:}
{\small
The aim of this paper is to firmly establish subjective fiducial inference as a rival to the more
conventional schools of statistical inference, and to show that Fisher's intuition concerning the
importance of the fiducial argument was correct. In this regard, methodology outlined in an earlier
paper is modified, enhanced and extended to deal with general inferential problems in which various
parameters are unknown.
As a key part of what is put forward, the joint fiducial distribution of all the parameters of a
given model is determined on the basis of the full conditional fiducial distributions of these
parameters by using an analytical approach or a Gibbs sampling method, the latter of which does not
require these conditional distributions to be compatible.
Although the resulting theory is classified as being `subjective', this is essentially due to the
argument that all probability statements made about fixed but unknown parameters must be inherently
subjective.
In particular, it is systematically argued that, in general, there is no need to place a great
emphasis on the difference between the fiducial probabilities derived by using this theory of
inference and objective probabilities.
Some important examples of the application of this theory are presented.}
\\[3ex]
\pagebreak

\section{Introduction}

R.\ A.\ Fisher is one of the greatest, if not the greatest, statistician that has ever lived.
Many of his contributions to statistical theory were considered to be revolutionary, but one
concept that he developed and discussed at length, namely the fiducial argument, has as yet failed
to gain many advocates. It is clear from his writings on the subject (see for example Fisher~1930,
1935, 1956) that he regarded this concept as being the foundation of a school of inference to
rival the other two main schools of inference that still flourish to this day, that is
Neyman\hspace{0.05em}-Pearson and Bayesian inference.
However, since Fisher's death in 1962, few have attempted to develop the fiducial argument into a
separate school of inference, notable exceptions being the work of D.\ A.\ S.\ Fraser on structural
inference, see Fraser~(1966, 1972), and the theory contained in Wilkinson~(1977).

There has been more activity, on the other hand, in attempting to use the fiducial argument to
support other schools of inference. In particular, it has been used to support
Dempster\hspace{0.03em}-Shafer theory, as is apparent in the original discussion of this theory,
e.g.\ see Dempster~(1968) and Shafer~(1976), and the Neyman\hspace{0.05em}-Pearson school of
inference by means of both generalized fiducial inference, see Hannig~(2009) and Hannig et
al.~(2016), and confidence distribution theory, see Xie and Singh~(2013). Also, it can be viewed as
supporting the theories of Dempster\hspace{0.03em}-Shafer and Neyman\hspace{0.05em}-Pearson
simultaneously as part of the recently developed theory of inferential models, see Martin and
Liu~(2015).

Nevertheless, for those who respect the intellect and intuition of Fisher, it may be disappointing
to see that one of his most cherished theories has been reduced to only a subsidiary role in
theories that have quite a distinct aim from what he had in mind.
Furthermore, the same people may be surprised by the fact that, even to this day, it is often
argued that fiducial inference is so closely related to Bayesian inference that if, as in many
cases, the fiducial distribution is equal to the posterior distribution for some choice of the
prior distribution, then the two theories are indistinguishable, see Lecoutre and
Poitevineau~(2014), Liu and Martin~(2015) and the work of many others.
In summary, it could be rather colourfully said that, in recent years, fiducial inference has been
like the wreck of a vintage car, which finds itself parked in a backstreet, sprayed with graffiti
by youths who do not appreciate its uniqueness and inner beauty, and robbed by opportunists for
spare parts to use in vehicles considered to be more commercially viable.

Following on from Bowater~(2017b), the aim of the present paper is to attempt to address this sorry
state of affairs. In particular, the theory of subjective fiducial inference put forward in
Bowater~(2017b) will be modified, enhanced and extended to deal with general inferential problems
in which various parameters are unknown.
As in this earlier paper, the type of inferences made by this theory will be regarded as being
`subjective'. However, the reason for this will be essentially attributed to the argument that all
probability statements made about fixed but unknown parameters must be inherently subjective,
rather than it being attributed to a particular need to emphasize how different fiducial
probabilities derived using this approach to inference are from objective probabilities.

Let us briefly outline the structure of the paper.
A motivation is given in the following section for the theory that will be developed, while the
concept of probability that this theory relies upon is explained in Section~\ref{sec7}.
After putting forward various important concepts, definitions and methods in Section~\ref{sec2},
the resulting approach to inference is applied to various examples in Section~\ref{sec12}.
The specification of the conceptual framework is then completed in Section~\ref{sec11}. The final
two sections of the paper clarify the merits of the theory and discuss some open issues.

\vspace{3ex}
\section{Motivation}
\label{sec9}

The need for an alternative theory of inference is motivated by the inadequacies of established
theories in addressing the issue of how to make inferences on the basis of data when nothing or
little was known about the parameters of interest before the data were observed. Here we will
briefly review the inadequacies of two such theories with regard to how they tackle this issue,
namely objective Bayesian inference and frequentist post-data inference.

Objective Bayesian inference is a form of Bayesian inference that is based on prior distributions
that have the property that the information contained within them has in some way been minimised,
either explicitly or implicitly, compared to the information that is expected to be contained in
the data. Advocates of this type of inference would argue that it offers a collection of methods
that attempt to standardise the way in which a prior distribution can be formed that, through the
posterior distribution, allows the data to `speak for themselves'. However, objective Bayesian
inference faces the following severe criticisms:

\vspace{2ex}
\noindent
1) If an objective Bayesian analysis is to avoid that inferences are dependent on the way the
sampling model is parameterised, which was a drawback of the classical Bayesian methods based on
the principle of insufficient reason that were proposed by both T.\ Bayes and P.\ S.\ Laplace, then
the choice of the prior distribution will need to depend on the sampling model, as is the case in
the methods that can be found, for example, in Jeffreys~(1961), Kass and Wasserman~(1996) and
Ghosh~(2011). However, as highlighted by many (see for example Seidenfeld~1979 and Lindley~1997),
this completely breaks the logic of Bayesian theory, since our state of knowledge about a parameter
will depend on how we intend to go about collecting more information about the parameter.

\vspace{2ex}
\noindent
2) Prior distributions derived using objective Bayesian methods are very often improper.
Such prior distributions break the standard rules of probability without having any special
permission for doing so, and are therefore purely mathematical creations that have no direct
real-world meaning.

\vspace{2ex}
\noindent
3) Even taking into account the principle of stable estimation (see Edwards, Lindman and
Savage~1963) and even when the sample size is large, posterior distributions will be generally very
sensitive to differences in the prior distributions that are derived by different objective
Bayesian methods.
Therefore, it is vital that there is a consensus on which objective Bayesian method is the best one
to use, but such a consensus does not exist. Some methods even lead to different prior
distributions depending on what is the parameter of main interest, e.g.\ the method outlined in
Bernardo~(1979) and Berger and Bernardo~(1992).

\vspace{2ex}
Let us now consider frequentist post-data inference, which we will assume refers to the type of
conditional frequentist approaches to inference outlined in Goutis and Casella~(1995).
The motivation for these methods clearly stems from the difficulties that arose from Fisher's
attempt to justify the fiducial argument in terms of frequentist probability. The reason that he
chose to do this would seem to have come from a desire to place the fiducial argument on an
objective footing and, in terms of quantifying uncertainty, objectivity for Fisher meant
frequentist probability.
To be more specific, the crux of Fisher's line of reasoning (as it is presented in Fisher~1956) was
that a fiducial interval for a parameter can only be valid in a frequentist sense if the sample
space contains no subsets that are \emph{recognisable} with respect to that interval.
However, Buehler and Feddersen~(1963) showed that even in one of the simplest and most common
problems of inference, that of making inferences about the mean $\mu$ of a normal distribution
when its variance $\sigma^2$ is unknown, recognisable subsets exist with respect to the standard
fiducial interval for $\mu$. To this day, the ubiquity of recognisable subsets in proposed
solutions to problems of inference represents a major obstacle to the further development of the
theory of frequentist post-data inference.

\vspace{2.5ex}
\section{A note about probability}
\label{sec7}

Let us begin by establishing the concept of probability that will underlie the theory that is about
to be developed.
We observe that it would be difficult to argue that subjective probability is a meaningless
concept. The fact that a meteorological expert can say that \linebreak he believes the probability
of rain tomorrow is 0.3 and others find this information useful would seem to satisfactorily refute
such an argument. Some would argue that subjective probability is the only concept of probability
that is required and, moreover, that making distinctions between different types of subjective
probability is essentially pointless.
Such a view is common amongst advocates of the Bayesian paradigm, see for example,
de Finetti~(1974, 1975) and Savage~(1954).

Perhaps the most standard position to take on whether probabilities are of different types is to
contend that two types of probability exist, one being subjective probability based on some kind of
elicitation method, and the other being frequentist probability based on calculating the long run
proportion of times a repeatable experiment produces a given outcome.
In this viewpoint, it would appear that the value that is assigned as the probability of any given
event is not sufficient to fully define the probability concerned, since we also need to know
whether the probability is subjective or frequentist.
Also, the prevalence of this viewpoint naturally gives importance to the
Bayesian\hspace{0.05em}-\hspace{0.03em}frequentist controversy, which arises due to the fact that,
according to specific criteria, inferences made using the Bayesian approach often conflict with
inferences made using the Neyman\hspace{0.05em}-Pearson (frequentist) approach.

This paper will rely on the definition of probability that was originally presented in
Bowater~(2017a) under the name of type B probability, and which was subsequently extended to
formally incorporate probability distributions in Bowater~(2017b).
Under this definition, probability comprises of two components, namely a probability value, which
is the sole recognised component in conventional definitions of probability, and the strength
assigned to this probability value.
Therefore, probabilities can be big and weak, small and strong, big and strong, small and weak,
etc. The strength component allows an ordered classification of probability types, and therefore is
more sophisticated than the standard dichotomous system of classifying probabilities as simply
being subjective or frequentist. For a full definition and explanation of this concept of
probability, the reader is referred to the two aforementioned papers. Nevertheless, to summarise
how this definition of probability can be used to determine a probability value and its strength
for a single event, a modified version of an example that appears in Bowater~(2017a) will now be
presented.

Let us suppose that an individual wishes to determine his probability for the event of a first-term
US president being re-elected in three years' time, which will be referred to as the event $A$.
From the earlier papers, it can be seen that we must first decide upon a reference set of events
$R = \{R_1,R_2,\ldots,R_m\}$. Taking into account the likely precision by which he may be able to
determine his probability value for the event $A$, let us imagine that the individual decides that
the events $R_i$ correspond to each of the outcomes of drawing a ball from an urn of 20 distinctly
labelled balls. With the event $R(\lambda)$ defined by substituting $m=20$ into the general
definition of this event, i.e.\
\vspace{1.25ex}
\begin{equation}
\label{equ23}
R(\lambda) = \left\{
\begin{array}{ll}
R_1 \cup R_2 \cup \cdots \cup R_{\lambda m} & \hspace{0.5em}\mbox{if $\lambda \in
\Lambda \cup \{ 1\}$}\\[1ex]
\emptyset & \hspace{0.5em}\mbox{if $\lambda=0$}
\end{array}
\right.
\vspace{1.25ex}
\end{equation}
where $\Lambda = \{1/m,\hspace{0.1em} 2/m,\ldots, (m-1)/m\}$, his probability value for the event
$A$ is then defined as being the unique value of $\lambda \in \{0, 0.05, 0.1, \ldots, 1\}$ that
maximises the similarity $S(A, R(\lambda))$, i.e.\ the similarity between his conviction that the
event $A$ will occur and his conviction that the event $R(\lambda)$ will occur.
Let us assume that this value is 0.7. Therefore, it is being assumed that the individual is capable
of asserting \pagebreak that, in his opinion, the similarities $S(A,R(0.65))$ and $S(A,R(0.75))$
are less than the similarity $S(A,R(0.7))$, which seems a reasonable assumption to make.

Now let us consider an event associated with spinning what is known as a probability wheel (see
Spetzler and Stael von Holstein~1975) which consists essentially of a rotatable disc with a fixed
pointer in its centre. Assuming that the area of the disc is divided into a red sector and a blue
sector, let the event of interest be the event of the pointer coming to rest in the red sector,
which will be referred to as the event $B$. If the proportion of the area of the disc that is red
is 0.7, then using the definition of probability being considered, it would not be at all
surprising if, with respect to the aforementioned reference set $R$, the individual assigned a
probability value of 0.7 to the event $B$.

However, although it will be assumed that a probability value of 0.7 will be assigned to both the
events $A$ and $B$, the strength that is associated with this probability value when it is assigned
to event $A$ is likely to be different from when it is assigned to event $B$, even under the
assumption that the probability value for the event $A$ has been determined as precisely as
possible.
In particular, it is likely that the similarity $S(A,R(0.7))$ will be considered to be
substantially less than the similarity $S(B,R(0.7))$, which is equivalent to asserting that the
probability of 0.7 is a much weaker probability for the event $A$ than for the event $B$.

The reason for this should be fairly evident, since the nature of the uncertainty about whether
event $A$ will occur is clearly different from the nature of the uncertainty about both whether
event $R(0.7)$ will occur and whether event $B$ will occur.
More specifically, the factors that can influence whether or not event $A$ will occur are likely to
be considered vague and difficult to weigh up, while events $R(0.7)$ and $B$, on the other hand,
are the outcomes of two standard types of physical experiment.

The idea that a probability comprises of both a probability value and its strength is supported by
the need to explain the expression of ambiguity aversion in decision making, which is an issue that
has been long debated in microeconomic theory (see Ellsberg~1961, Gilboa and Schmeidler~1989 and
Alary, Gollier and Treich~2013).
In this regard, the definition of probability under discussion has been used to undermine the
independence axiom upon which the foundations of Bayesian theory depend, since it facilitates a
rational explanation of the paradox associated with Ellsberg's three colour example (see
Ellsberg~1961 for this example and Bowater~2017a for the explanation).
Therefore, this counters the popular argument that any measure of the uncertainty of an event that
is not solely the probability value assigned to the event must be invalid due to the measure not
being compatible with Bayesian theory.

The concept of strength can be applied not just to individual probabilities, but also to entire
probability density functions in the sense that, under additional criteria, one density function
can be classified as being weaker or stronger than another density function.
Loosely speaking, a probability density $f_{X}(x)$ is defined as being stronger than another
density $g_{\hspace{0.1em}Y}(y)$ at the level of resolution $\alpha$, if probabilities equal to
$\alpha$ derived by integrating $f_{X}(x)$ over subspaces of $x$ are considered to be at least as
strong as, and sometimes stronger than, probabilities equal to $\alpha$ derived by integrating
$g_{\hspace{0.1em}Y}(y)$ over subspaces of $y$.
A more detailed definition of this property can be found in Bowater~(2017b).

Although in the sections that immediately follow attention will be focused on the determination of
probability values and densities rather than on the determination of their strengths, this latter
issue needs to be borne in mind.
We will explicitly return to the task of completing probability definitions by assigning strengths
to probability values and densities in Section~\ref{sec11}.

\vspace{3ex}
\section{Subjective fiducial inference}
\label{sec2}

Let us now present in detail, what will be called, the theory of subjective fiducial inference.

\vspace{3ex}
\subsection{Sampling model}

In general, it will be assumed that a sampling model that depends on one or various unknown
parameters $\theta = \{\theta_i: i=1,2,\ldots,k\}$ generates the data $x=\{x_i: i=1,2,\ldots,n\}$.
Let the joint density of the data given the true values of the parameters $\theta$ be denoted as
$g(x\,|\,\theta)$.

\vspace{3ex}
\subsection{Univariate case}
\label{sec8}

For the moment, we will assume that the only unknown parameter in the model is
$\theta_j$\hspace{0.05em}, either because there are no other parameters in the model, or because
the true values of the parameters $\theta_{-j} = \{ \theta_1,\ldots,\theta_{j-1},\theta_{j+1},
\ldots,\theta_k \}$ are known.

\vspace{3ex}
\noindent
{\bf Definition 1: Fiducial statistics}

\vspace{1ex}
\noindent
Given this assumption, a fiducial statistic $Q(x)$ will be defined as being a one-dimensional
sufficient statistic for $\theta_j$ if such a statistic exists, otherwise it may be assumed to be
any one-to-one function of a unique maximum likelihood estimator of $\theta_j$.

\vspace{3ex}
\noindent
{\bf Assumption 1: Data generating algorithm}

\vspace{1ex}
\noindent
Independent of the way in which the data were actually generated, it will be assumed that the data
set $x$ was generated by the following algorithm:

\vspace{2ex}
\noindent
1) Generate a value $\gamma$ for a continuous one-dimensional random variable $\Gamma$, which has a
probability density function $f_{\Gamma}(\gamma)$ that does not depend on the parameter
$\theta_j$\hspace{0.05em}.

\vspace{2ex}
\noindent
2) Determine a value $q(x)$ for a fiducial statistic $Q(x)$ by setting $\Gamma$ equal to $\gamma$
and $Q(x)$ equal to $q(x)$ in the specific form of the definition of the distribution of $Q(x)$
that is given by:
\begin{equation}
\label{equ1}
Q(x)=\varphi(\Gamma, \theta_j)
\end{equation}
where the function $\varphi(\Gamma, \theta_j)$ is defined so that it satisfies the following
conditions:

\vspace{2ex}
\noindent
{\bf Assumption~1.1: Conditions on the function} $\varphi(\Gamma, \theta_j)$

\vspace{1ex}
\noindent
a) The distribution of $Q(x)$ as defined by equation~(\ref{equ1}) is equal to what it would have
been if $Q(x)$ had been determined on the basis of the data set $x$.\\
b) The only random variable upon which $\varphi(\Gamma, \theta_j)$ depends is the variable
$\Gamma$.\\
c) Let $G = \{ \gamma : f_{\Gamma}(\gamma) > 0 \}$, and let $H$ be the set of all possible values
of $\theta_j$ as specified before any information about the data $x$ has been obtained.
If it is assumed that a value for $Q(x)$ has been generated, but both its corresponding value
$\gamma$ for the variable $\Gamma$ and the parameter $\theta_j$ are unknown, then substituting
$Q(x)$ in equation~(\ref{equ1}) by whatever value is taken by $Q(x)$ would imply that this equation
would define an injective mapping from the set $G$ to the set $H$.

\vspace{2ex}
\noindent
3) Generate the data set $x$ from the sampling density $g(x\,|\,\theta_1,\theta_2, \ldots,
\theta_k)$ conditioned on the statistic $Q(x)$ being equal to its already generated value $q(x)$.

\vspace{3ex}
In the context of the above algorithm, the variable $\Gamma$ will be referred to as a primary
random variable (primary r.v.). However, if the above algorithm was rewritten so that the value
$\gamma$ of the variable $\Gamma$ was generated by setting it equal to a deterministic function of
an already generated value for $Q(x)$ and the parameter $\theta_j$\hspace{0.05em}, then $\Gamma$
would not be a primary r.v. In relation to Neyman\hspace{0.05em}-Pearson theory, a primary r.v.\
could be classified as a type of pivot that is distinguished in terms of the way it is generated.

\vspace{3ex}
\noindent
{\bf Definition 2: Univariate subjective fiducial distributions}

\vspace{1ex}
\noindent
Given a value $q(x)$ for the fiducial statistic $Q(x)$, the subjective fiducial distribution of
\linebreak the parameter $\theta_j$ conditional on all other parameters $\theta_{-j}$ being known
is defined by set\-ting $Q(x)$ equal to $q(x)$ in equation~(\ref{equ1}), and then treating the
value $\theta_j$ in this equation as being a realisation of the random variable
$\Theta_j$\hspace{0.05em}, to give the expression:
\begin{equation}
\label{equ3}
q(x)=\varphi(\Gamma, \Theta_j)
\end{equation}
where $\Gamma$ has the density function $f_{\Gamma}(\gamma)$ defined in step~1 of the data
generating algorithm in Assumption~1.
As can easily be shown, this fiducial distribution for $\theta_j$ does not de\-pend on the choice
made for the statistic $Q(x)$.
Also, observe that condition (c) of Assumption~1.1 ensures that equation~(\ref{equ3}) specifies a
valid probability distribution for the parameter $\theta_j$.
The classical fiducial argument can be seen through the fact that the distribution of the primary
r.v.\ $\Gamma$ is the same both before and after the fiducial statistic $Q(x)$ is observed.

\vspace{3ex}
\subsection{Multivariate case}
\label{sec3}

We will now consider the case where all the parameters $\theta = \{ \theta_1, \theta_2, \ldots,
\theta_k\}$ in the sampling model are unknown.

For any given data set $x$, let us assume that the method outlined in the previous section allows
us to define the fiducial density of the parameter $\theta_j$ conditional on all other parameters
$\theta_{-j}$ for all values of $j$, i.e.\ $j=1,2,\ldots, k$. We will denote this set of full
conditional fiducial densities as:
\begin{equation}
\label{equ2}
f(\theta_j\,|\,\theta_{-j},x)\ \ \ \ \mbox{for $j=1,2,\ldots,k$}
\end{equation}
If these conditional densities determine a unique joint density for all the parameters $\theta$,
then this density will be defined as being the joint subjective fiducial density of these
parameters and will be denoted as $f(\theta\,|\,x)$.
However, the set of densities in equation~(\ref{equ2}) may not be consistent with any joint density
of the parameters concerned, i.e.\ these conditional densities may be incompatible among
themselves. On the other hand, if the conditional densities under discussion are indeed compatible
then, since, under a mild requirement, a joint density function is uniquely defined by its full
conditional densities, these densities will, in general, define a unique joint fiducial density for
the parameters~$\theta$.

Therefore, it would be helpful to know the cases in which the conditional densities in
equation~(\ref{equ2}) are compatible and when they are not, and if they are indeed incompatible,
whether and how the difficulty that this leads to can be addressed.
In this regard, we will propose two different strategies. The first strategy is to establish
whether the full conditional densities in question are compatible using analytical methods.
By contrast, the second strategy is to assume that these conditional densities are incompatible
even when they may not be, and use a computational method to try to find the joint density function
of all the parameters $\theta$ that has full conditional densities that most closely approximate
the densities in equation~(\ref{equ2}). We now will discuss each of these strategies in a bit more
detail.

\vspace{3ex}
\subsection{Verifying the compatibility of full conditional distributions}
\label{sec1}

Various analytical methods have been proposed for establishing the compatibility of full
conditional distribution functions in a general context, see for example Arnold and Press~(1989),
Arnold, Castillo and Sarabia~(2002) and Kuo and Wang~(2011). Nevertheless, these methods can
largely only be applied to cases where the variables over which these distribution functions are
defined can only take a finite number of different values, or where there are only two such
variables. There are, though, two methods of this type that at least potentially are more widely
applicable. Therefore, we now will take a look at these two methods.

The first method we will consider is a simple one.
In particular, we begin by proposing an analytical expression for the joint density function of the
set of parameters $\theta$, then we determine the full conditional density functions for this joint
density, and finally we see whether these conditional densities are equivalent to the full
conditional densities in equation~(\ref{equ2}). If this equivalence is achieved, then these latter
conditional densities clearly must be compatible. This method has the advantage that, in such
circumstances, it directly gives us, under a mild condition, an analytical expression for the
unique joint fiducial density of the parameters $\theta$, i.e.\ under this condition, it will be
the originally proposed joint density for these parameters.
Since this joint fiducial density may well be equal to a joint posterior density of the parameters
$\theta$ that, given some choice for the joint prior density of these parameters, is derived using
Bayes' theorem, a good proposal for this joint fiducial density, at least up to a normalising
constant, may often be found by multiplying the likelihood function by a convenient mathematical
choice for the joint prior density of the parameters $\theta$.

The second method that we will consider for verifying the compatibility of the set of full
conditional densities of interest depends on studying the behaviour of a Gibbs sampling algorithm
(Geman and Geman~1984, Gelfand and Smith~1990) that makes transitions on the basis of this set of
conditional densities.
In particular, let us define a single transition of this type of algorithm as being one that
results from randomly drawing a value (only once) from each of the full conditional densities in
equation~(\ref{equ2}) according to some given fixed order of these densities, which we will call a
fixed scanning order, replacing each time the previous value of the parameter concerned by the
value that is generated.
To clarify, it is being assumed that only the set of values for the parameters $\theta$ that are
obtained on completing a transition of this kind are recorded as being a newly generated sample,
i.e.\ the intermediate sets of parameter values that are used in the process of making such a
transition do not form part of the output of the algorithm.

On the basis of the results in Chen and Ip~(2015), it can be deduced that the full conditional
densities in equation~(\ref{equ2}) will be compatible if the Gibbs sampling algorithm just
described satisfies the following three conditions:

\vspace{1.5ex}
\noindent
1) It is positive recurrent for all possible fixed scanning orders. This condition ensures that the
sampling algorithm has at least one stationary distribution for any given fixed scanning order.

\vspace{1.5ex}
\noindent
2) It is irreducible and aperiodic for all possible fixed scanning orders. Together with
condition~1, this condition ensures that the sampling algorithm has a limiting distribution for any
given fixed scanning order.

\vspace{1.5ex}
\noindent
3) Given conditions~1 and~2 hold, the limiting density function of the sampling algorithm needs to
be the same over all possible fixed scanning orders.

\vspace{1.5ex}
\noindent
Moreover, when these conditions hold, the joint fiducial density function of the parameters
$\theta$ implied by the full conditional densities under discussion will be the unique limiting
density function of these parameters referred to in condition~3.
The sufficiency of the conditions~1 to~3 just listed for establishing the compatibility of any
given set of full conditional densities was proved for a special case in Chen and Ip~(2015), which
is a proof that can be easily extended to the more general case that is currently of interest.

In the context of the full conditional densities of the parameters $\theta$ being the full
conditional fiducial densities in equation~(\ref{equ2}), let us briefly comment on how easy it is
likely to be, in practice, to establish whether or not these conditional densities satisfy each of
the three conditions in question.
First, it would not be expected to be that difficult, in this context, to determine whether or not
condition~1 is satisfied, since a failure of this condition to hold would only be expected to occur
in very pathological cases.
Also, the fulfilment of condition~2 will usually be easy to verify through an inspection of the
full conditional densities concerned.
On the other hand, in the context of interest, it will usually be impossible to determine whether
or not condition~3 is satisfied.
Despite this substantial drawback, we will nevertheless consider again the strategy that has just
been outlined in the next subsection.

\vspace{3ex}
\subsection{Finding compatible approximations to incompatible full conditionals}
\label{sec14}

In any given situation where it is not easy to establish whether or not the full conditional
fiducial densities in equation~(\ref{equ2}) are compatible, let us imagine that we make the
pessimistic assumption that they are in fact incompatible.
Nevertheless, even though these conditional densities could be incompatible, they could be 
reasonably assumed to represent the best information that is available for constructing a joint
density function for the parameters $\theta$ that most accurately represents what is known about
these parameters after the data have been observed, i.e.\ constructing, what could be referred to
as, the most suitable joint fiducial density for these parameters.
Therefore, it would seem appropriate to try to find the joint density of the parameters concerned
that has full conditional densities that most closely approximate those given in
equation~(\ref{equ2}).

Various methods have been proposed for doing this outside of the context of the full conditional
densities being fiducial densities and when the random variables concerned can only take a finite
number of different values, which means of course that we need to refer to the probability mass
rather than density functions of these variables, see for example Arnold, Castillo and
Sarabia~(2002), Chen, Ip and Wang~(2011), Chen and Ip~(2015) and Kuo, Song and Jiang~(2017).
Similar to what was done earlier, here we will again focus attention on a more widely applicable
method, in particular the method that simply consists in making the assumption that the joint
density of the parameters $\theta$ that most closely corresponds to the set of full conditional
densities in equation~(\ref{equ2}) is equal to the limiting density function of a Gibbs sampling
algorithm that is based on these conditional densities with some given fixed or random scanning
order of the parameters concerned.
This approach relates to more specific methods for addressing the general problem of interest that
were discussed in Chen, Ip and Wang~(2011) and Mur\'e~(2019).
To clarify, a transition of the Gibbs sampler in question under a random scanning order will be
defined as being one that results from generating a value from one of the conditional densities in
equation~(\ref{equ2}) that is chosen at random, with the probability of any given density
$f(\theta_j\,|\,\theta_{-j},x)$ being selected being set equal to some given value $a_j$, where of
course $\sum_{i=1}^{k} a_i = 1$, and then treating the generated value as the updated value of the
parameter concerned.

To measure how close the full conditional densities of the limiting density function of the general
type of Gibbs sampler under discussion are to the full conditional densities in
equation~(\ref{equ2}), we can use a variation on the line of reasoning that underlies the second
method for verifying the compatibility of full conditional densities that was outlined in the last
subsection.
In particular, assuming that condition~1 (positive recurrence condition) and condition~2
(irreducibility and aperiodicity condition) of this method are satisfied, it would appear to be
useful (with reference to condition~3 of this method) to analyse how the limiting density function
of the Gibbs sampler being considered varies over a reasonable number of very distinct fixed
scanning orders of the sampler.
If within such an analysis, the variation of this limiting density with respect to the scanning
order of the parameters $\theta$ can be classified as small, negligible or undetectable, then this
should give \linebreak us reassurance that the full conditional densities in equation~(\ref{equ2})
are, respectively according to such classifications, close, very close or at least very close, to
the full conditional densities of the limiting density of a Gibbs sampler of the type that is of
main interest, i.e.\ a Gibbs sampler that is based on any given fixed or random scanning order of
the parameters concerned.

In trying to choose the scanning order of this type of Gibbs sampler such that it has a limiting
density function that corresponds to a set of full conditional densities that most accurately
approximate the densities in equation~(\ref{equ2}), a good general choice would arguably be the
random scanning order of the parameters $\theta$ that was defined earlier with the selection
probability of any given parameter, i.e.\ the probability $a_j$, being set equal to
$1/k$\hspace{0.05em} for all $j$, which is what we will call a uniform random scanning order.
In a loosely similar context, Mur\'e~(2019) not only recommends but provides some analytical
results to support the use of a uniform random scanning order within such a sampling procedure to
address the issue in question.

However, it can be easily shown that independent of whether or not the set of full conditional
densities in equation~(\ref{equ2}) are compatible, the last full conditional density in this set
that is sampled from in completing a given fixed scanning order will be one of the full conditional
densities of the limiting density function of the type of Gibbs sampler being discussed that uses
such a fixed scanning order.
This therefore provides a reason for perhaps deciding, in certain applications, that the limiting
density of a Gibbs sampler of the type of interest most satisfactorily corresponds to the full
conditional densities in equation~(\ref{equ2}) when a given fixed rather than a uniform random
scanning order of the parameters $\theta$ is used.

As with all Gibbs samplers it is important to verify in implementing any of the aforementioned
strategies that the sampler concerned has converged to its limiting density function within the
restricted number of transitions of the sampler that can be observed in practice.
To do this we can make use of standard methods for analysing the convergence of Monte Carlo Markov
chains described in, for example, Gelman and Rubin~(1992), Cowles and Carlin~(1996) and Brooks and
Gelman~(1998).
However, the use of such convergence diagnostics may be considered to be slightly more important in
the case of present interest in which the full conditional densities on which the Gibbs sampler is
based could be incompatible, since, compared to the case where these densities are known to be
compatible, there is likely to be, in practice, a little more concern that the Gibbs sampler may
not actually have a limiting density function, even though in reality the genuine risk of this may
still be extremely low.

A notable advantage of the general method for finding a suitable joint fiducial density for the
parameters $\theta$ that has just been outlined is that it can directly achieve what is often the
main goal of a standard application of the Gibbs sampler, namely that of obtaining good
approximations to the expected values of functions of the parameters of a model over the post-data
or posterior density for these parameters that is of interest, i.e.\ expected values of the
following type:
\vspace{1ex}
\[
\mbox{E}[h(\theta)\,|\,x] = \int_{\mbox{\footnotesize{$\mathbb{R}^k$}}} h(\theta)p(\theta\,|\,x)
d\theta
\vspace{1ex}
\]
where $p(\theta\,|\,x)$ is a given post-data density function of the parameters $\theta$, which in
the current context would of course be the joint fiducial density of these parameters, i.e.\
$f(\theta\,|\,x)$, while $h(\theta)$ is any given function of the parameters concerned.
To be more specific, this kind of expected value may, of course, be approximated using the Monte
Carlo estimator:
\vspace{2ex}
\[
\frac{1}{m-b}\hspace{0.1em} \sum_{i\hspace{0.05em}=\hspace{0.05em}b\hspace{0.05em}+1}^{m}
h(\theta_1^{(i)},\theta_2^{(i)},\ldots,\theta_k^{(i)})
\vspace{2ex}
\]
where $\theta_1^{(i)},\theta_2^{(i)},\ldots,\theta_k^{(i)}$ is the $i$th sample of parameter values
among the $m$ samples generated by the sampler in total, and $b$ is the number of initial samples
that are classified as belonging to the burn-in phase of the sampler.

Finally, it is worth noting that when the sampling model has only two parameters, i.e.\ $k=2$, it
is easy to show that the limiting marginal densities of the two parameters concerned that are
produced by running a Gibbs sampler that is based on incompatible full conditional densities of
these two parameters are not affected by the scanning order of the sampler over the types of fixed
and random scanning orders of the sampler that were defined earlier.
In the current context, this property may be of some convenience if the aim is to only determine
the marginal fiducial densities of the parameters of a model that has only two parameters.
It is, though, a property that does not generally hold when there are more than two parameters.

\vspace{3ex}
\section{Applications to multivariate cases}
\label{sec12}
In this section, applications will be presented of the methodology detailed in the preceding
sections to cases where the sampling model contains more than one unknown parameter.

\vspace{3ex}
\subsection{Inference about the mean and variance of a normal distribution}
\label{sec5}

To begin with, let us consider the standard problem of making inferences about the mean $\mu$ of a
normal distribution, when its variance $\sigma^2$ is unknown, on the basis of a sample $x$ of size
$n$, i.e.\ $x=\{x_1,x_2,\ldots,x_n\}$, drawn from the distribution concerned.
Although a solution to this problem using subjective fiducial inference was put forward in
Bowater~(2017b), it should become clear from what is about to be presented and what is discussed in
Section~\ref{sec11} that the approach outlined in the current paper provides a more elegant way in
which this problem can be resolved.

In the first instance, let us assume that $\sigma^2$ is known. Under this assumption, a sufficient
statistic for $\mu$ would be the sample mean $\bar{x}$, which therefore will be treated as being
the fiducial statistic $Q(x)$ in this particular case.
Defining the primary r.v.\ $\Gamma$ as having a standard normal distribution implies that, given a
value for $\sigma^2$, equation~(\ref{equ1}) can be expressed as:
\vspace{1.5ex}
\[
\bar{x}=\varphi(\Gamma,\mu)=\mu+(\sigma/\sqrt{n}\hspace{0.1em})\hspace{0.05em}\Gamma
\vspace{1.5ex}
\]
meaning that, according to Definition~2 of Section~\ref{sec8}, the fiducial distribution of $\mu$
is defined by:
\vspace{1ex}
\begin{equation}
\label{equ5}
\mu\, |\, \sigma^2, x \sim \mbox{N}(\bar{x}, \sigma^2 / n)
\vspace{1ex}
\end{equation}
which is the standard fiducial distribution of $\mu$ for this problem conditional on $\sigma^2$
being known.
On the other hand, if $\mu$ was known, a sufficient statistic \pagebreak for $\sigma^2$ would be
the variance estimator:
\vspace{1.25ex}
\begin{equation}
\label{equ17}
\bm\hat{\sigma}^2 = (1/n) \textstyle{\sum_{i=1}^{n} (x_i-\mu)^2}
\vspace{1.25ex}
\end{equation}
which will therefore be treated as being the statistic $Q(x)$ in this particular case. Defining
$\Gamma$ as having a $\chi^2$ distribution with $n$ degrees of freedom implies that, given a value
for $\mu$, equation~(\ref{equ1}) can be expressed as:
\[
\bm\hat{\sigma}^2 = \varphi(\Gamma,\sigma^2)=(\sigma^2/n)\Gamma
\]
meaning that the fiducial distribution of $\sigma^2$ is defined by:
\begin{equation}
\label{equ6}
\sigma^2\, |\, \mu, x \sim \mbox{Scale-inv-$\chi^2$} (n, \bm\hat{\sigma}^2)
\end{equation}
i.e.\ it is a scaled inverse $\chi^2$ distribution with $n$ degrees of freedom and scaling
parameter equal to $\bm\hat{\sigma}^2$, which is the standard fiducial distribution of $\sigma^2$
for this problem conditional on $\mu$ being known.

To verify that the full conditional distributions of $\mu$ and $\sigma^2$ in equations~(\ref{equ5})
and~(\ref{equ6}) determine a joint distribution for $\mu$ and $\sigma^2$, we can use the simple
analytical method outlined in the opening part of Section~\ref{sec1}.
In particular, the full conditional distributions of the joint proper posterior distribution of
$\mu$ and $\sigma^2$ that corresponds to choosing the prior density of $\mu$ and $\sigma^2$ to be
the improper density $p(\mu,\sigma^2) \propto 1/\sigma^2$ are identical to the full conditionals in
equations~(\ref{equ5}) and~(\ref{equ6}).
Therefore, the conditional distributions in these equations are compatible, and it is clear that
the joint fiducial distribution for $\mu$ and $\sigma^2$ that they define must be unique. More
specifically, by integrating over this joint distribution with respect to $\sigma^2$, it can be
established that the marginal fiducial distribution of $\mu$ is defined by:
\[
\mu\,|\,x \sim \mbox{Non-standardised}\ t_{n-1} (\bar{x}, s/\sqrt{n}\hspace{0.1em})
\]
where $s$ is the sample standard deviation, i.e.\ it is the familiar \pagebreak non-standardised
Student $t$ distribution with $n-1$ degrees of freedom, location parameter equal to $\bar{x}$ and
scaling parameter equal to $s/\sqrt{n}$, which indeed is the standard marginal fiducial
distribution of $\mu$ for the problem of interest.

\vspace{3ex}
\subsection{Inference about both parameters of a Pareto distribution}

Let us now consider the problem of making inferences about both the shape parameter $\alpha$ and
the scale parameter $\beta$ of a Pareto distribution on the basis of a sample $x$ of size $n$ from
the density function concerned, i.e.\ the function:
\vspace{1.5ex}
\[
g(y\,|\,\alpha,\beta) = \left\{
\begin{array}{ll}
\alpha \beta^{\alpha} y^{-(\alpha+1)} \ & \mbox{if $y \geq \beta$}\\[1ex]
0 & \mbox{otherwise}
\end{array}
\right.
\vspace{1.5ex}
\]
where $y$ is any given value in the sample $x$.

If $\beta$ was known, a sufficient statistic for $\alpha$ would be $\sum_{i=1}^{n} \log x_i$, which
therefore will be treated as being the fiducial statistic $Q(x)$ in this particular case.
Defining the primary r.v.\ $\Gamma$ as having a $\mbox{Gamma}(n,1)$ distribution, i.e.\ a gamma
distribution with shape parameter equal to $n$ and rate parameter equal to 1, implies that, given a
value for $\beta$, equation~(\ref{equ1}) can be expressed as:
\vspace{0.5ex}
\[
\textstyle{\sum_{i=1}^{n}} \log x_i = \varphi(\Gamma,\alpha)=(\Gamma/\alpha)+n\log \beta
\vspace{0.5ex}
\]
meaning that, according to Definition~2, the fiducial distribution of $\alpha$ is defined by:
\vspace{0.5ex}
\begin{equation}
\label{equ7}
\alpha\, |\, \beta, x \sim \mbox{Gamma}\hspace{0.1em} (\hspace{0.1em}
\mbox{scale} = n,\hspace{0.2em} \mbox{rate} = \textstyle{\sum_{i=1}^{n}} (\log x_i -
\log \beta )\hspace{0.1em} )
\vspace{0.5ex}
\end{equation}
On the other hand, if $\alpha$ was known, a sufficient statistic for $\beta$ would be the minimum
value of the sample, i.e.\ $\min(x)$, which will therefore be treated as being the statistic $Q(x)$
in this particular case.
Defining $\Gamma$ as having an exponential \pagebreak distribution with rate parameter equal to 1,
implies that, given a value for $\alpha$, equation~(\ref{equ1}) can be expressed as:
\[
\min(x) = \varphi(\Gamma,\beta) = \exp((\Gamma/n\alpha)+\log\beta)
\]
meaning that the fiducial density of $\beta$ is defined by:
\vspace{1.75ex}
\begin{equation}
\label{equ8}
f(\beta \,|\, \alpha, x) = \left\{
\begin{array}{ll}
(n\alpha/\beta) \exp \{ -n\alpha (\log(\min(x)) - \log\beta) \} \ \ & \mbox{if $0 < \beta \leq
\min(x)$}\\[1ex]
0 & \mbox{otherwise}
\end{array}
\right.
\vspace{1.75ex}
\end{equation}

The full conditional distributions of the joint proper posterior distribution of $\alpha$ and
$\beta$ that corresponds to choosing the prior density of $\alpha$ and $\beta$ to be the improper
density $p(\alpha,\beta) \propto (\alpha\beta)^{-1}$ are identical to the full conditional
distributions of $\alpha$ and $\beta$ in equations~(\ref{equ7}) and~(\ref{equ8}).
Therefore, the conditional distributions in these equations are compatible and it is clear that the
joint fiducial density of $\alpha$ and $\beta$ that they determine must be unique.
In particular, from what was just mentioned, it can be deduced that this joint fiducial density
must be defined by:
\vspace{1.5ex}
\[
f(\alpha, \beta \,|\, x) = \left\{
\begin{array}{ll}
{\tt C}\hspace{0.1em} \alpha^{n-1} \beta^{n\alpha-1} \prod_{i=1}^{n} x_i^{-(\alpha+1)} \ \ &
\mbox{if $0 < \beta \leq \min(x)$ and $\alpha>0$}\\[1ex]
0 & \mbox{otherwise}
\end{array}
\right.
\vspace{1.5ex}
\]
where ${\tt C}$ is a normalising constant.

\vspace{3ex}
\subsection{Inference about all parameters of a normal quadratic regression model}

To show how the approach outlined in Section~\ref{sec2} can be applied to make inferences about the
parameters of normal polynomial regression models, let us consider the example of applying this
approach to the problem of making inferences about all the parameters $\beta_0$, $\beta_1$,
$\beta_2$ and $\sigma^2$ of a normal quadratic regression model, i.e.\
\vspace{0.5ex}
\[
y_{i} = \beta_0 + \beta_1 x_i + \beta_2 x_i^2 + \varepsilon_i\hspace{0.2em},\ \ \ \mbox{where
$\varepsilon_i \sim \mbox{N}(0,\sigma^2)$}
\vspace{0.5ex}
\]
on the basis of a sample $y_{+}=\{(x_i, y_i):i=1,2,\ldots,n\}$ \pagebreak from the model concerned.

Sufficient statistics for each of the parameters $\beta_0$, $\beta_1$, $\beta_2$ and $\sigma^2$
conditional on all parameters except the parameter itself being known are respectively:
\vspace{0.25ex}
\[
\textstyle{
\sum_{i=1}^{n} y_i,\ \hspace{0.2em}\sum_{i=1}^{n} x_i y_i,\
\hspace{0.2em}\sum_{i=1}^{n} x_i^2 y_i\ \hspace{0.3em} \mbox{and}\
\hspace{0.2em}\sum_{i=1}^{n} (y_i - \beta_0 - \beta_1 x_i - \beta_2 x_i^2)^2}
\]
which would therefore be suitable choices for the fiducial statistic $Q(y_{+})$ for the particular
cases in question.
By applying the methodology described in Section~\ref{sec2} under the assumption that each of these
statistics in turn is treated as being the statistic $Q(y_{+})$, it can be shown that the full
conditional fiducial distributions for this problem are defined by:
\vspace{1ex}
\begin{equation}
\label{equ9}
\beta_0\,|\,\beta_1,\beta_2,\sigma^2,y_{+} \sim \mbox{N}\hspace{-0.1em} \left(\hspace{0.1em}
\mbox{$\sum_{i=1}^{n}$}\hspace{0.1em} y_i/n - \beta_{1}\hspace{0.1em}
\mbox{$\sum_{i=1}^{n}$}\hspace{0.1em} x_i/n - \beta_{2}\hspace{0.1em}
\mbox{$\sum_{i=1}^{n}$}\hspace{0.1em} x_i^{2}/n,\hspace{0.25em} \sigma^2/n \hspace{0.1em}\right)
\vspace{2.5ex}
\end{equation}
\begin{equation}
\label{equ10}
\beta_1\,|\,\beta_0,\beta_2,\sigma^2,y_{+} \sim \mbox{N}\hspace{-0.1em} \left(
\frac{\sum_{i=1}^{n} x_{i}y_{i} - \beta_{0}\sum_{i=1}^{n} x_i - \beta_{2}\sum_{i=1}^{n} x_i^{3}}
{\sum_{i=1}^{n} x_{i}^2},\hspace{0.25em} \frac{\sigma^2}{\sum_{i=1}^{n} x_{i}^2} \right)
\vspace{2.5ex}
\end{equation}
\begin{equation}
\label{equ11}
\beta_2\,|\,\beta_0,\beta_1,\sigma^2,y_{+} \sim \mbox{N}\hspace{-0.1em} \left(
\frac{\sum_{i=1}^{n} x_{i}^{2}y_{i} - \beta_{0}\sum_{i=1}^{n} x_i^{2} -
\beta_{1} \sum_{i=1}^{n} x_i^{3}}{\sum_{i=1}^{n} x_{i}^4},\hspace{0.25em}
\frac{\sigma^2}{\sum_{i=1}^{n} x_{i}^4} \right)
\vspace{2.5ex}
\end{equation}
\begin{equation}
\label{equ12}
\sigma^2\, |\, \beta_0, \beta_1, \beta_2,y_{+} \sim \mbox{Scale-inv-$\chi^2$}\hspace{-0.1em}
\left(\hspace{0.1em}n,\hspace{0.15em} (1/n)\hspace{0.1em} \mbox{$\sum_{i=1}^{n}$}\hspace{0.1em}
(y_i - \beta_0 - \beta_1 x_i - \beta_2 x_i^2)^2 \hspace{0.1em}\right)
\vspace{2ex}
\end{equation}

The full conditional distributions of the joint proper posterior distribution of $\beta_0$,
$\beta_1$, $\beta_2$ and $\sigma^2$ that corresponds to choosing the joint prior density of these
parameters to be the improper density $p(\beta_0,\beta_1,\beta_2,\sigma^2) \propto 1/\sigma^2$ are
identical to the full conditionals in equations~(\ref{equ9}) to~(\ref{equ12}).
Therefore, the conditional distributions in these equations are compatible, and from the
information just given, it can deduced that they determine a unique joint fiducial density for
$\beta_0$, $\beta_1$, $\beta_2$ and $\sigma^2$ that is defined by:
\[
f(\beta_0, \beta_1, \beta_2, \sigma^2 \,|\, y_{+})\, \propto\, \sigma^{-(n+2)} \exp \left\{
\hspace{0.1em} - (1/2\sigma^2)\hspace{0.1em} \mbox{$\sum_{i=1}^{n}$}\hspace{0.1em}
( y_i-\beta_{0} - \beta_{1}x_i -\beta_{2}x_{i}^2)^2 \hspace{0.1em}\right\}
\]

\vspace{3ex}
\subsection{Inference about both parameters of a gamma distribution}
\label{sec4}

We will now consider the problem of making inferences about both the shape parameter $\alpha$ and
the rate parameter $\beta$ of a gamma distribution on the basis of a sample $x$ of size $n$ from
the density function concerned, i.e.\ the function:
\vspace{1ex}
\[
g(y\,|\,\alpha,\beta) =
\frac{\beta^{\alpha} y^{\alpha-1} \exp(-y\beta)}{\mathtt{G}(\alpha)}\ \ \ \ \mbox{if $y \geq 0$
and zero otherwise}
\vspace{1ex}
\]
where $y$ is any given value in the sample $x$ and the function $\mathtt{G}(\alpha)$ is the gamma
function evaluated at $\alpha$.

If $\alpha$ was known, a sufficient statistic for $\beta$ would be $\sum_{i=1}^{n} x_i$, which
therefore will be treated as being the fiducial statistic $Q(x)$ in this particular case.
Defining the primary r.v.\ $\Gamma$ as having a $\mbox{Gamma}(n\alpha,1)$ distribution, i.e.\ a
gamma distribution with $\mbox{shape} = n\alpha$ and $\mbox{rate}=1$, implies that, given a value
for $\alpha$, equation~(\ref{equ1}) can be expressed as:
\[
\textstyle{\sum_{i=1}^{n}} x_i = \varphi(\Gamma,\beta)=\Gamma/\beta
\]
meaning that the fiducial distribution of $\beta$ is defined by:
\begin{equation}
\label{equ13}
\beta\, |\, \alpha, x \sim \mbox{Gamma}\hspace{0.1em} (\hspace{0.1em} \mbox{shape} = n\alpha,
\hspace{0.2em} \mbox{rate} = \textstyle{\sum_{i=1}^{n}}\hspace{0.1em} x_i \hspace{0.1em} )
\vspace{0.25ex}
\end{equation}

On the other hand, if $\beta$ was known, a sufficient statistic for $\alpha$ would be
$\sum_{i=1}^{n} \log x_i$, which will therefore be treated as being the statistic $Q(x)$ in this
particular case.
However, defining the primary r.v.\ $\Gamma$ and the required function $\varphi(\Gamma,\alpha)$ in
equation~(\ref{equ1}) is not straightforward in this case. This is due to the cumulative density
function of $\sum_{i=1}^{n} \log x_i$ not being mathematically very tractable. A rudimentary way of
approximating the distribution of $\sum \log x_i$ is to use the central limit theorem. This is the
approximation method that will be adopted here.

In this regard, it is relevant to note that the mean and variance of $\sum \log x_i$ can be shown
to be $n(\psi(\alpha)-\log\beta)$ and $n\psiprime(\alpha)$ respectively, \pagebreak where
$\psi(\alpha)$ and $\psiprime(\alpha)$ are, respectively, the digamma and trigamma functions
evaluated at $\alpha$.
Therefore, assuming that $\beta$ is known and that $\sum \log x_i$ is approximately normally
distributed, equation~(\ref{equ1}) can be approximated by:
\vspace{0.5ex}
\begin{equation}
\label{equ14}
\textstyle{\sum_{i=1}^{n}} \log x_i = \varphi(\Gamma,\alpha) = n(\psi(\alpha)-\log\beta) +
\Gamma \sqrt{n\psiprime(\alpha)}
\vspace{0.5ex}
\end{equation}
where $\Gamma$ is defined as having a $\mbox{N}(0,1)$ distribution.
If $n$ is sufficiently large then, given a fixed value of $\sum \log x_i$, this equation defines an
injective mapping from subsets of values $\gamma$ for the variable $\Gamma$ to the space of
$\alpha$ except when these subsets of $\gamma$ values contain extremely positive or negative
numbers.
Therefore, the function $\varphi(\Gamma,\alpha)$ in this equation approximately satisfies
condition~(c) of Assumption~1.1, and as a result, under Definition~2, this equation approximately
defines the fiducial distribution of $\alpha$ conditional on $\beta$ being known.

To illustrate this example, Figure~1 shows some results from running a Gibbs sampler on the basis
of the full conditional fiducial distributions of the parameters $\beta$ and $\alpha$ defined by
equations~(\ref{equ13}) and~(\ref{equ14}) respectively, with a uniform random scanning order of
these two parameters, as such a scanning order was defined in Section~\ref{sec14}.
In particular, the histograms in Figures~1(c) and~1(d) represent the distributions of the values of
$\alpha$ and $\beta$, respectively, over a single run of five million samples of these parameters
generated by the Gibbs sampler after a preceding run of two thousand samples, which were classified
as belonging to its burn-in phase, had been discarded.
The data set $x$ was a random sample of $n=30$ values drawn from a gamma distribution with
$\alpha=2$ and $\beta=0.5$.
Also, Figures~1(a) and~1(b) show the progression of the sampler over the first two thousand
`cycles' of its main (rather than its burn-in) sampling run in terms of the parameters $\alpha$ and
$\beta$ respectively, i.e.\ its progression in generating the first 2,000 samples of this run.

\begin{figure}[!t]
\begin{center}
\noindent
\makebox[\textwidth]{\includegraphics[width=7in]{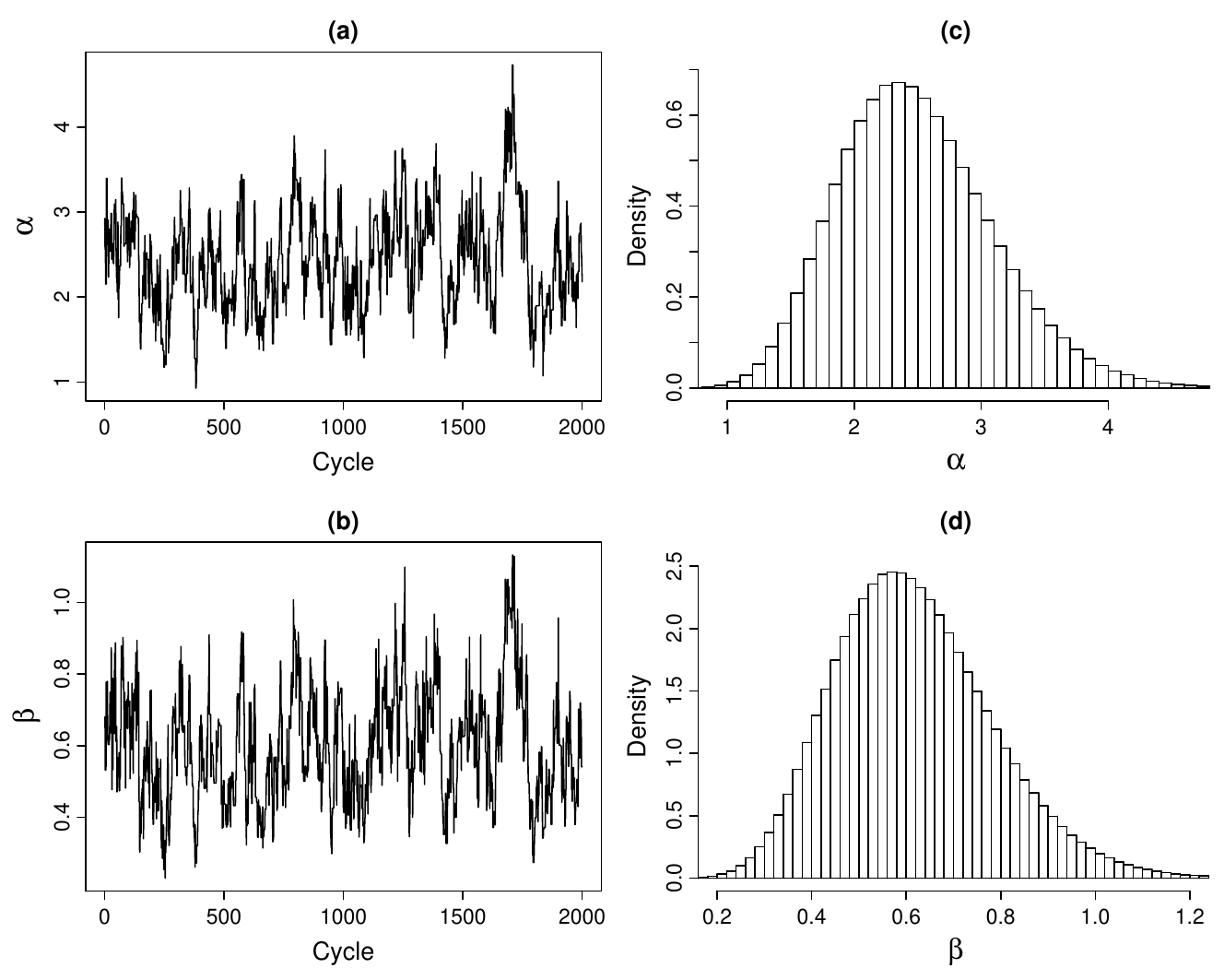}}
\caption{\small{Gibbs sampling based on the full conditional fiducial distributions of the
parameters of a gamma distribution}}
\end{center}
\end{figure}

To generate each value of $\alpha$ from the full conditional fiducial distribution of $\alpha$
defined by equation~(\ref{equ14}), a random value $\gamma$ of the variable $\Gamma$ was first
\pagebreak drawn from a $\mbox{N}(0,1)$ distribution truncated to lie in the interval $[-5,5]$, and
then equation~(\ref{equ14}) was numerically solved to find the corresponding value of the parameter
$\alpha$.
Truncating the distribution of $\Gamma$ in the way just described meant that there was always an
injective mapping from the space of possible values for $\gamma$ (i.e.\ the values $-5$ to 5) to
the space of $\alpha$, i.e.\ condition~(c) of Assumption~1.1 was always satisfied.

In addition to this analysis, the Gibbs sampler was run various times from different starting
points, and a careful study of the output of these runs using appropriate diagnostics provided no
evidence to suggest that the sampler does not have a limiting distribution, and showed, at the same
time, that it would appear to generally converge quickly to this distribution.
Furthermore, the Gibbs sampling algorithm was run separately with each of the two possible fixed
scanning orders of the parameters, i.e.\ the one in which $\alpha$ is updated first and then
$\beta$ is updated, and the one that has the reverse order, in accordance with how a single
transition of such an algorithm was defined in Section~\ref{sec1}, i.e.\ single transitions of the
algorithm incorporated updates of both parameters.
In doing this, no statistically significant difference was found between the samples of parameter
values aggregated over the runs of the sampler in using each of these two scanning orders after
excluding the burn-in phase of the sampler, e.g.\ between the two sample correlations of $\alpha$
and $\beta$, even when the runs concerned were long.
Taking into account what was discussed in Section~\ref{sec14}, this implies that the full
conditional distributions of the limiting distribution of the original Gibbs sampler, i.e.\ the one
with a uniform random scanning order, should be, at the very least, close approximations to the
full conditional distributions on which the sampler is based, i.e.\ the conditional fiducial
distributions defined by equations~(\ref{equ13}) and~(\ref{equ14}).

\vspace{3ex}
\subsection{Inference about both parameters of a beta distribution}
\label{sec6}

The next problem we will consider is that of making inferences about both parameters $\alpha$ and
$\beta$ of a beta distribution on the basis of a sample $x$ of size $n$ from the density function
concerned, i.e.\ the function:
\vspace{0.5ex}
\[
g(y\,|\,\alpha,\beta) = \frac{y^{\alpha-1} (1-y)^{\beta-1}}{\mathtt{B}(\alpha,\beta)}\ \ \ \
\mbox{if $0 \leq y \leq 1$ and zero otherwise}
\vspace{0.5ex}
\]
where $y$ is any given value in the sample $x$ and the function $\mathtt{B}(\alpha,\beta)$ is the
beta function evaluated at $\alpha$ and $\beta$.

If $\beta$ was known, a sufficient statistic for $\alpha$ would be $\sum_{i=1}^{n} \log x_i$, which
therefore will be treated as being the fiducial statistic $Q(x)$ in this particular case. However,
the cumulative density function of $Q(x)$ in this case is again not mathematically very tractable.
For this reason, similar to what was done in the previous example, the central limit theorem will
be used to approximate the distribution of the statistic $Q(x)$.

On making the assumption, which we will therefore be choosing to make, that the statistic
$\sum_{i=1}^{n} \log x_i$ is approximately normally distributed, it follows that, given a value for
$\beta$, equation~(\ref{equ1}) can be approximated by:
\begin{eqnarray}
\hspace{-1em}\textstyle{\sum_{i=1}^{n}} \log x_i = \varphi(\Gamma,\alpha)\hspace{-0.5em} & = &
\hspace{-0.5em}\mbox{E}[Q(x)] + \mbox{Var}[Q(x)]^{0.5}\hspace{0.1em} \Gamma \nonumber\\
& = & \hspace{-0.5em}n(\psi(\alpha)-\psi(\alpha+\beta))
+ n^{0.5}(\psiprime(\alpha)-\psiprime(\alpha+\beta))^{0.5}\hspace{0.1em} \Gamma \label{equ15}
\vspace{1.25ex}
\end{eqnarray}
where $\Gamma \sim \mbox{N}(0,1)$ and the functions $\psi(b)$ and $\psiprime(b)$ are again the
digamma and trigamma functions evaluated at any given value $b$.
If $n$ is sufficiently large, then the function $\varphi(\Gamma,\alpha)$ in this equation satisfies
condition (c) of Assumption~1.1 under the restriction that values of $\gamma$ that are extremely
positive or negative are excluded from the set $G$. As a result, under Definition~2, this equation
approximately defines the fiducial distribution of $\alpha$ conditional on $\beta$ being known.

On the other hand, if $\alpha$ was known, a sufficient statistic for $\beta$ would be
$\sum_{i=1}^{n} \log (1-x_i)$, which will therefore be treated as being the statistic $Q(x)$ in
this particular case.
Since the cumulative density function of $Q(x)$ in this case is again not mathematically very
tractable, the distribution of $Q(x)$ will be once more approximated on the basis of the central
limit theorem.
Having made the assumption, which we will therefore be choosing to make, that $\sum_{i=1}^{n}
\log (1-x_i)$ is approximately normally distributed, it follows that, given a value for $\alpha$,
equation~(\ref{equ1}) can be approximated by:
\begin{eqnarray}
\hspace{-2em}\textstyle{\sum_{i=1}^{n} \log (1-x_i)} = \varphi(\Gamma,\beta) \hspace{-0.5em} & =
& \hspace{-0.5em}\mbox{E}[Q(x)] + \mbox{Var}[Q(x)]^{0.5}\hspace{0.1em} \Gamma \nonumber\\
& = & \hspace{-0.5em}n(\psi(\beta)-\psi(\alpha+\beta))
+ n^{0.5}(\psiprime(\beta)-\psiprime(\alpha+\beta))^{0.5}\hspace{0.1em} \Gamma \label{equ16}
\end{eqnarray}
where again $\Gamma \sim \mbox{N}(0,1)$.
Similar to the previous case, the function $\varphi(\Gamma,\beta)$ in this equa-{\linebreak}tion
approximately satisfies condition~(c) of Assumption~1.1 and, as a result, under Definition~2, this
equation approximately defines the fiducial distribution of $\beta$ conditional on $\alpha$ being
known.

To illustrate this example, Figure~2 shows some results from running a Gibbs sampler on the basis
of the full conditional fiducial distributions of the parameters $\alpha$ and $\beta$ defined by
equations~(\ref{equ15}) and~(\ref{equ16}) respectively, with a uniform random scanning order of
these two parameters.
In particular, the histograms in Figures~2(c) and~2(d) represent the distributions of the values of
$\alpha$ and $\beta$, respectively, over a single run of five million samples of these parameters
generated by the Gibbs sampler after a preceding run of two thousand samples were discarded due to
these samples being classified as belonging to its burn-in phase.
The data set $x$ was a random sample of $n=50$ values drawn from a beta distribution with
$\alpha=8$ and $\beta=3$.
Also, Figures~2(a) and~2(b) show the progression of the sampler over the first two thousand
samples/cycles of its main sampling run in terms of the parameters $\alpha$ and $\beta$
respectively.
Each random value of $\alpha$ and $\beta$ was generated from the full conditional fiducial
distributions for $\alpha$ and $\beta$ under discussion using the same numerical method that was
used in Section~\ref{sec4} to generate random values from the fiducial distribution defined by
equation~(\ref{equ14}), meaning therefore that, as a first step, a random value of the variable
$\Gamma$ was drawn from a $\mbox{N}(0,1)$ distribution truncated to lie in the interval $[-5,5]$.
The reason for truncating the distribution of $\Gamma$ in this way was the same as given in
Section~\ref{sec4}.

\begin{figure}[!t]
\begin{center}
\noindent
\makebox[\textwidth]{\includegraphics[width=7in]{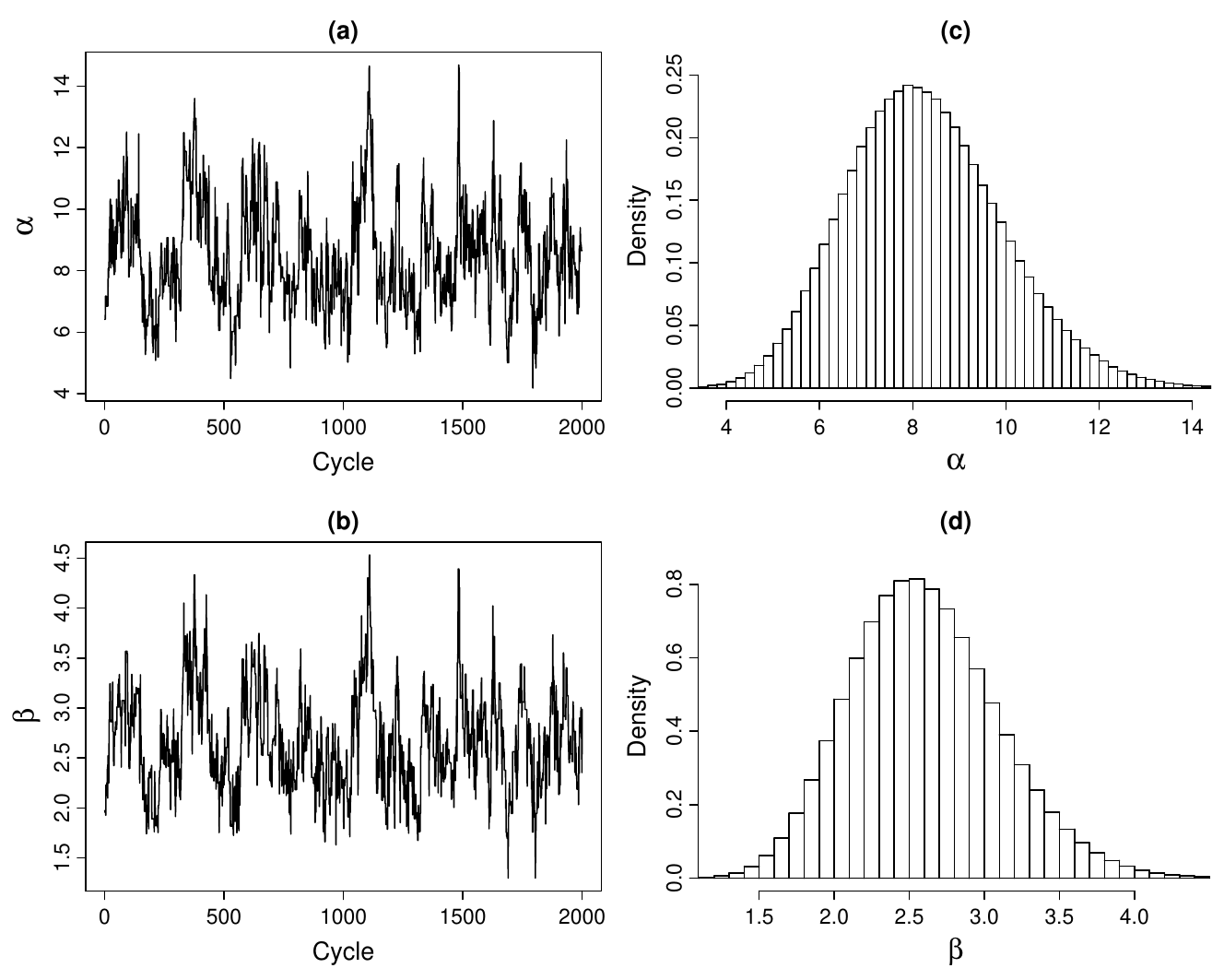}}
\caption{\small{Gibbs sampling based on the full conditional fiducial distributions of the
parameters of a beta distribution}}
\end{center}
\end{figure}

Supplementary to this analysis, the Gibbs sampler was run various times from different starting
points, and there was no suggestion from using appropriate diagnostics that the sampler does not
have a limiting distribution.
Furthermore, after excluding the burn-in phase of the sampler, no statistically significant
difference was found between the samples of parameter values aggregated over the runs of the
sampler in using each of the two fixed scanning orders of the parameters $\alpha$ and $\beta$ that
are possible, with a single transition of the sampler defined in the same way as in the example
outlined in the previous section, even when the runs concerned were long.
Therefore, taking into account what was discussed in Section~\ref{sec14}, the full conditional
distributions of the limiting distribution of the original random-scan Gibbs sampler should be, at
the very least, close approximations to the full conditional distributions on which the sampler is
based, i.e.\ the conditional fiducial distributions defined by equations~(\ref{equ15})
and~(\ref{equ16}).

\vspace{3ex}
\subsection{The Behrens\hspace{0.05em}-Fisher problem}

As a prelude to examining the general problem of making inferences about all the parameters of a
bivariate normal distribution based on a data set consisting of realisations of the two random
variables $X$ and $Y$ described by this distribution, let us consider the special case of this
problem in which the covariance of $X$ and $Y$ is assumed to be zero. The data set will be assumed
to comprise of a sample $x$ of $n_x$ realisations of the variable $X$ and a sample $y$ of $n_y$
realisations of the variable $Y$. Let $\mu_x$ and $\mu_y$ denote the means of the variables $X$ and
$Y$ respectively, and let the variances of $X$ and $Y$ be denoted by $\sigma_x^2$ and $\sigma_y^2$
respectively.

On the basis of the results presented in Section~\ref{sec5}, it is clear that the full conditional
fiducial distributions for this problem are defined by first substituting $\mu_x$ for $\mu$,
$\sigma_x^2$ for $\sigma^2$ and $n_x$ for $n$ into equations~(\ref{equ5}), (\ref{equ17})
and~(\ref{equ6}) to obtain the conditional fiducial distributions
$f(\mu_x\,|\,\mu_y,\sigma_x^2,\sigma_y^2,x,y)$ and $f(\sigma_x^2\,|\,\mu_x,\mu_y,\sigma_y^2,x,y)$,
and then by substituting $\mu_y$ for $\mu$, $\sigma_y^2$ for $\sigma^2$, $n_y$ for $n$ and the
sample $y$ for the sample $x$ into the same equations to get the conditional fiducial distributions
$f(\mu_y\,|\,\mu_x,\sigma_x^2,\sigma_y^2,x,y)$ and $f(\sigma_y^2\,|\,\mu_x,\mu_y,\sigma_x^2,x,y)$.

The full conditional distributions of the joint proper posterior distribution of $\mu_x$, $\mu_y$,
$\sigma_x^2$ and $\sigma_y^2$ that corresponds to choosing the prior density of $\mu_x$, $\mu_y$,
$\sigma_x^2$ and $\sigma_y^2$ to be the improper density\hspace{0.1em} $p(\mu_x, \mu_y,\sigma_x^2,
\sigma_y^2)$ $\propto$ $1/\sigma_x^2\sigma_y^2$\hspace{0.1em} are identical to the full conditional
fiducial distributions that have just been defined.
Therefore, these latter conditional distributions are compatible, and from the information just
given, it can deduced that on the basis of the unique joint fiducial distribution of $\mu_x$,
$\mu_y$, $\sigma_x^2$ and $\sigma_y^2$ that they determine, the marginal fiducial distribution of
$\mu_x-\mu_y$ is defined by:
\vspace{1ex}
\begin{equation}
\label{equ18}
\mu_x - \mu_y = \bar{x} - \bar{y} + B \sqrt{(s_x^2/n_x) + (s_y^2/n_y)}
\pagebreak
\end{equation}
where $s_x^2$ and $s_y^2$ are the observed variances of the samples $x$ and $y$ respectively, and
$B$ is a random variable that has a Behrens\hspace{0.05em}-Fisher distribution with degrees of
freedom $n_x-1$ and $n_y-1$ (order irrelevant), and with angle parameter equal to\hspace{0.05em}
$\tan^{-1} ((s_x \sqrt{n_y}) / (s_y \sqrt{n_x}))$. The distribution for $\mu_x - \mu_y$ defined by
equation~(\ref{equ18}) is the fiducial distribution for $\mu_x - \mu_y$ that was advocated by
R.\ A.\ Fisher for this problem.

\vspace{3ex}
\subsection{Inference about all parameters of a bivariate normal distribution}
\label{sec10}

The final problem we will consider in this section is the problem of making inferences about all
five parameters of a bivariate normal distribution, i.e.\ the means $\mu_x$ and $\mu_y$ and the
variances $\sigma^2_x$ and $\sigma^2_y$, respectively, of the two random variables concerned $X$
and $Y$, and the correlation $\rho$ of $X$ and $Y$, on the basis of a sample from this type of
distribution, i.e.\ the sample $\mathtt{z}=\{(x_i,y_i) : i=1,2,\ldots,n\}$, where $x_i$ and $y_i$
are the $i$th realisations of $X$ and $Y$ respectively.

If all parameters except $\mu_x$ were known, a sufficient statistic for $\mu_x$ would be:
\vspace{1ex}
\[
\sum^{n}_{i=1} x_i - \rho \left( \frac{\sigma_x}{\sigma_y} \right) \sum^{n}_{i=1} y_i
\vspace{1ex}
\]
which therefore will be treated as being the fiducial statistic $Q(\mathtt{z})$ in this particular
case. Defining the primary r.v.\ $\Gamma$ as having a $\mbox{N}(0,1)$ distribution, implies that,
given the values of all parameters except $\mu_x$, equation~(\ref{equ1}) can be expressed as:
\vspace{2ex}
\[
\sum^{n}_{i=1} x_i - \rho \left( \frac{\sigma_x}{\sigma_y} \right) \sum^{n}_{i=1} y_i =
\varphi(\Gamma,\mu_x) = n\mu_x- n\rho \left( \frac{\sigma_x}{\sigma_y} \right) \mu_y +
(n\sigma_x^2(1-\rho^2))^{0.5}\hspace{0.15em} \Gamma
\vspace{2ex}
\]
meaning that the fiducial distribution of $\mu_x$ is defined by:
\vspace{1.5ex}
\begin{equation}
\label{equ21}
\mu_x\, |\,\mu_y,\sigma_x^2,\sigma_y^2,\rho, {\tt z} \sim \mbox{N} \left(\hspace{0.1em} \bar{x} +
\rho\hspace{-0.05em} \left(\frac{\sigma_x}{\sigma_y} \right)\hspace{-0.1em}
(\mu_y - \bar{y}),\hspace{0.4em} \frac{\sigma_x^2(1-\rho^2)}{n}\hspace{0.1em} \right)
\pagebreak
\end{equation}
Due to the symmetrical nature of the bivariate normal distribution, it should be clear that, if all
parameters except $\mu_y$ are known, the fiducial distribution of $\mu_y$ is defined by:
\vspace{2ex}
\begin{equation}
\label{equ22}
\mu_y\, |\,\mu_x,\sigma_x^2,\sigma_y^2,\rho, {\tt z} \sim \mbox{N} \left(\hspace{0.1em} \bar{y} +
\rho\hspace{-0.05em} \left(\frac{\sigma_y}{\sigma_x} \right)\hspace{-0.1em}
(\mu_x - \bar{x}),\hspace{0.4em} \frac{\sigma_y^2(1-\rho^2)}{n}\hspace{0.1em} \right)
\vspace{2.5ex}
\end{equation}

By contrast, if all parameters except $\sigma_x^2$ are known, then no single sufficient statistic
for $\sigma^2_x$ exists, and therefore, in agreement with Definition~1 of Section~\ref{sec8}, we
will define the fiducial statistic $Q(\mathtt{z})$ to be the unique maximum likelihood estimator of
$\sigma_x^2$ given all other parameters are known.
This estimator is the value $\bm\hat{\sigma}_x^2$ that solves the following quadratic equation:
\vspace{2ex}
\begin{equation}
\label{equ25}
n(1-\rho^2)\bm\hat{\sigma}_x^2 + \rho \left(\frac{\bm\hat{\sigma}_x}{\sigma_y}
\right) \sum_{i=1}^{n} x'_i y'_i = 0
\vspace{2ex}
\end{equation}
where $x'_i = x_i - \mu_x$ and $y'_i = y_i - \mu_y$.

Now, it is well known that a maximum likelihood estimator of a parameter is usually asymptotically
normally distributed with mean equal to the true value of the parameter and variance equal to the
inverse of the Fisher information with respect to that parameter (i.e.\ the Fisher information
obtained via differentiating the logarithm of the likelihood function with respect to that
parameter). Since it can be shown that the Fisher information of the likelihood function in this
example with respect to $\sigma_x$, assuming all other parameters are known, is given by:
\vspace{1ex}
\[
{\cal I}(\sigma_x) = \frac{n(2-\rho^2)}{\sigma^2_x (1 - \rho^{\hspace{0.02em}2})}\ ,
\vspace{1.5ex}
\]
equation~(\ref{equ1}) can therefore be approximated, in the case of interest, by:
\vspace{2ex}
\[
\bm\hat{\sigma}_x = \sqrt{\varphi(\Gamma,\sigma^2_x)} = \sigma_x + \Gamma \sigma_x \left(
\frac{(1-\rho^2)}{n(2-\rho^{\hspace{0.02em}2})} \right)^{\hspace{-0.1em}0.5}
\vspace{2ex}
\]
where $\bm\hat{\sigma}_x$ is the maximum likelihood estimator of $\sigma_x$ defined by \pagebreak
equation~(\ref{equ25}) and $\Gamma$ is defined as having a $\mbox{N}(0,1)$ distribution.
Solving this equation for $\sigma^2_x$ leads to the follow\-ing approximate definition of the
fiducial distribution for $\sigma^2_x$ given all other parameters are known:
\vspace{3ex}
\begin{equation}
\label{equ19}
\sigma^2_x = \bm\hat{\sigma}^2_x \left(\Gamma \left(
\frac{(1-\rho^2)}{n(2-\rho^{\hspace{0.02em}2})} \right)^{\hspace{-0.1em}0.5} + 1
\right)^{\hspace{-0.1em}-2}
\vspace{3ex}
\end{equation}
Again due to the symmetrical nature of the bivariate normal distribution, the fiducial distribution
of $\sigma^2_y$ conditional on all parameters except $\sigma^2_y$ being known is approximately
defined in the same way, except that $\sigma^2_x$ and $\bm\hat{\sigma}^2_x$ in this expression are
substituted by $\sigma^2_y$ and $\bm\hat{\sigma}^2_y$ respectively, where $\bm\hat{\sigma}^2_y$ is
the maximum likelihood estimator of $\sigma^2_y$ given the values of all other parameters.

In looking now at the case where all parameters except the correlation coefficient $\rho$ are
known, we must again note that no single sufficient statistic exists for the unknown parameter of
interest, i.e.\ the parameter $\rho$ in this case, and therefore, similar to the case just
discussed, we will define the fiducial statistic $Q(\mathtt{z})$ to be the unique maximum
likelihood estimator of $\rho$ given all other parameters are known.
This estimator is the value $\bm\hat{\rho}$ that solves the following cubic equation:
\vspace{1.5ex}
\[
-n\bm\hat{\rho}^{\hspace{0.02em}3} + \left( \frac{\sum_{i=1}^n x'_i y'_i}{\sigma_x \sigma_y}
\right)\hspace{-0.1em} \bm\hat{\rho}^{\hspace{0.02em}2} + \left( n -
\frac{\sum_{i=1}^n (x'_i)^2}{\sigma_x^2} - \frac{\sum_{i=1}^n (y'_i)^2}{\sigma_y^2} \right)
\hspace{-0.1em} \bm\hat{\rho}\hspace{0.05em} + \hspace{0.03em}
\frac{\sum_{i=1}^n x'_i y'_i}{\sigma_x \sigma_y} = 0
\vspace{1.5ex}
\]
The distribution of this estimator will be approximated in the same way as the distributions of the
estimators $\bm\hat{\sigma}_x^2$ and $\bm\hat{\sigma}_y^2$ were approximated.
In particular, since it can be shown that the Fisher information of the likelihood function with
respect to $\rho$, assuming all other parameters are known, is given by:
\vspace{0.5ex}
\[
{\cal I}(\rho) = \frac{n(1+\rho^2)}{(1-\rho^{\hspace{0.02em}2})^2}\ ,
\vspace{0.5ex}
\]
equation~(\ref{equ1}) can be approximated, in the \pagebreak case in question, by:
\begin{equation}
\label{equ20}
\bm\hat{\rho} = \varphi(\Gamma,\rho) = \rho + \frac{(1 - \rho^2)\Gamma}{\sqrt{n(1 +
\rho{\hspace{0.02em}^2})}}
\vspace{2ex}
\end{equation}
where again $\Gamma \sim \mbox{N}(0,1)$. Under Definition~2, this equation defines the fiducial
distribution for $\rho$ given all other parameters are known.
It can be shown that if a random value of $\Gamma$ is substituted into equation~(\ref{equ20}), then
the value of $\rho$ that solves this equation will be unique. This value of $\rho$ will be, of
course, a random value of $\rho$ from the fiducial distribution of interest.

To illustrate this example, Figure~3 shows some results from running a Gibbs sampler on the basis
of the full conditional fiducial distributions of the parameters $\mu_x$, $\mu_y$, $\sigma_x^2$,
$\sigma_y^2$ and $\rho$ that were defined in the equations just given, i.e.\
equations~(\ref{equ21}) and~(\ref{equ22}), equation~(\ref{equ19}) in terms of both $\sigma_x^2$ and
$\sigma_y^2$, and equation~(\ref{equ20}), with a uniform random scanning order of these five
parameters.
In particular, the histograms in Figures~3(f) to~3(i) represent the distributions of the values of
$\mu_x$, $\mu_y$, $\sigma_x$, $\sigma_y$ and $\rho$, respectively, over a single run of five
million samples of these parameters generated by the Gibbs sampler after allowing for its burn-in
phase by discarding a preceding run of five thousand samples.
The data set $\mathtt{z}$ was a random sample of $n=200$ values drawn from a bivariate normal
distribution with $\mu_x=\mu_y=0$, $\sigma^2_x = \sigma^2_y = 1$ and $\rho=0.8$.
Also, Figures~3(a) to~3(e) show the progression of the sampler over the first five thousand
samples/cycles of its main sampling run in terms of each of the model parameters.

\begin{figure}[p]
\begin{center}
\noindent
\makebox[\textwidth]{\includegraphics[width=5.55in]{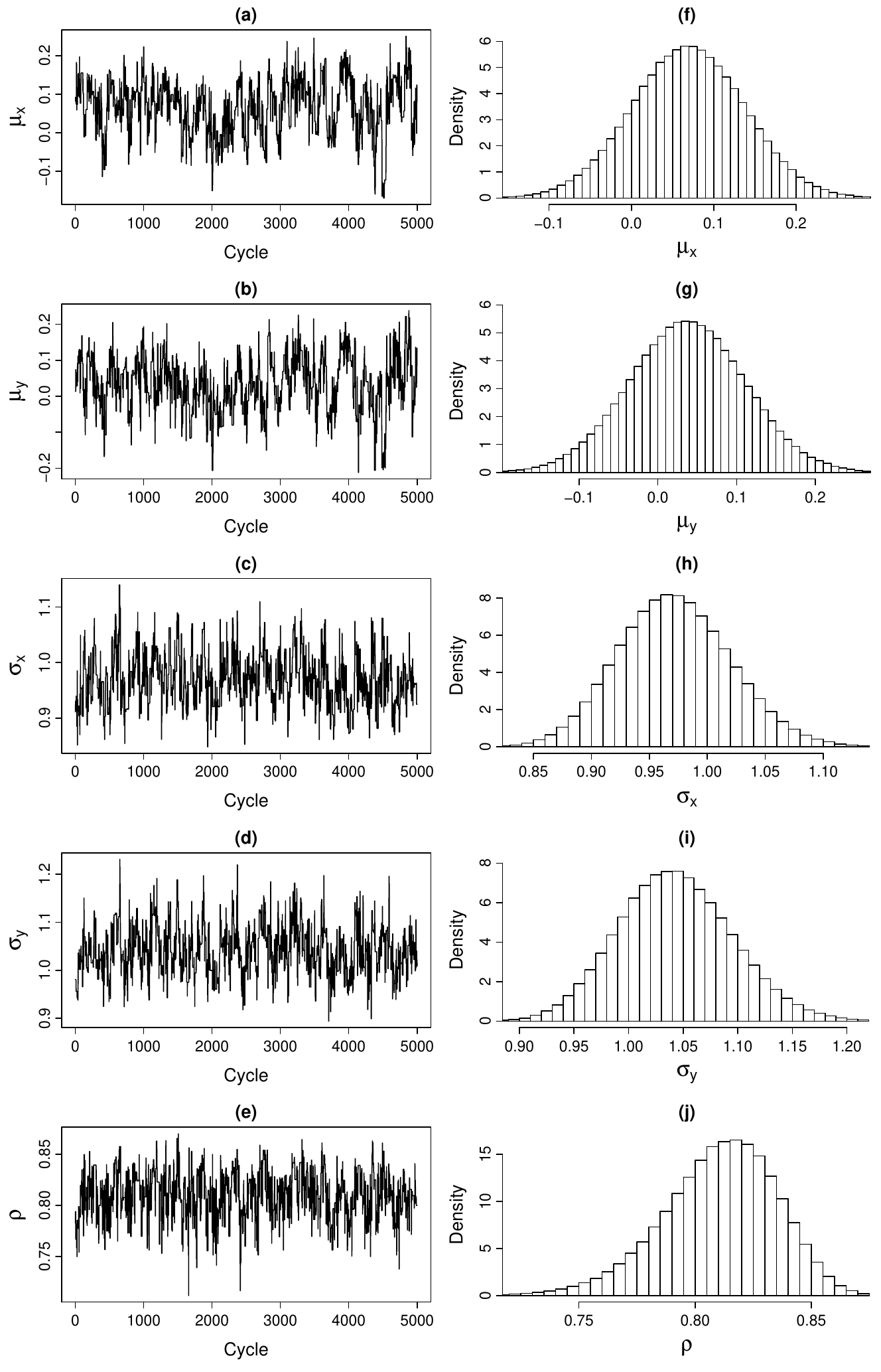}}
\vspace{-4ex}
\caption{\small{Gibbs sampling based on the full conditional fiducial distributions of the
parameters of a bivariate normal distribution}}
\end{center}
\end{figure}

Supplementary to this analysis, there was no suggestion from applying appropriate diagnostics to
multiple runs of the Gibbs sampler from different starting points that it did not have a limiting
distribution.
Furthermore, the Gibbs sampling algorithm was run separately with various very distinct fixed
scanning orders of the five model parameters $\mu_x$, $\mu_y$, $\sigma_x^2$, $\sigma_y^2$ and
$\rho$ in accordance with how a single transition of such an algorithm with a fixed scanning order
was defined in Section~\ref{sec1}.
In doing this, no statistically significant difference was found between the samples of parameter
\pagebreak values aggregated over moderately long runs of the sampler in using each of the scanning
orders concerned after excluding the burn-in phase of the sampler, e.g.\ between the various
correlation matrices of the parameters and between the various distributions of each individual
parameter.
However, over longer runs, important differences were found between these aggregated samples of
parameter values in using particular fixed scanning orders of the parameters. More specifically,
variations in the correlation of $\rho$ and $\sigma_x$ and the correlation of $\rho$ and $\sigma_y$
were statistically significant between the runs of samples that corresponded to these fixed
scanning orders.

Therefore, on the grounds of what was discussed in Section~\ref{sec14}, it can be concluded that
while the full conditional fiducial distributions on which the Gibbs sampling algorithm is based
are almost certainly incompatible, these conditional distributions should still be fairly close
approximations to the full conditional distributions of the limiting distribution of the original
random-scan Gibbs sampler, i.e.\ this latter joint distribution of the \linebreak five model
parameters $\mu_x$, $\mu_y$, $\sigma_x^2$, $\sigma_y^2$ and $\rho$ should be consistent to a
reasonably precise degree with the full conditional fiducial distributions of these parameters that
were directly specified in the present section.
In addition, it would not seem justifiable to regard this issue as being a major drawback of the
method that has been outlined, taking into account that three of these latter full conditional
densities were only derived in an approximate manner on the basis of a specific normality
assumption, i.e.\ the full conditional fiducial distributions of $\sigma_x^2$, $\sigma_y^2$ and
$\rho$.

\vspace{3ex}
\section{Determining the strengths of subjective fiducial probabilities}
\label{sec11}

We will now fulfil the undertaking made at the end of Section~\ref{sec7} and complete, in effect,
the definitions of the fiducial densities that have been derived in the preceding sections by
drawing some general conclusions about what should be the strengths that are assigned to the
probability values that are obtained by integrating over these densities.

Using again the terminology proposed in Bowater~(2017a), let the reference set $R$ be the
balls-in-an-urn reference set defined in Section~\ref{sec7} but with $m$ instead of 20 balls in the
urn, which we will assume are numbered from 1 to $m$.
Observe that if, within the methodology of this earlier paper, it is quite reasonably assumed that
the similarity $S$ between any given event $A$ and another event is maximised when the other event
is also $A$, then the event $R(\lambda)$ defined by equation~(\ref{equ23}) must have a probability
with respect to the set $R$ that equals $\lambda$.

To make a comparison with this event, let us consider the event of the primary r.v.\ $\Gamma$ being
less than $\gamma(\lambda)$, where $\gamma(\lambda)$ is defined by:
\vspace{1.5ex}
\[
\int_{-\infty}^{\gamma(\lambda)} f_\Gamma (\gamma) d\gamma = \lambda
\vspace{1.5ex}
\]
in which the density $f_\Gamma (\gamma)$ is defined as in Assumption~1 and the value $\lambda$ is a
member of the set $\Lambda$ defined immediately below equation~(\ref{equ23}).
The probability $\lambda$ that would be assigned both to this latter event
$\{\Gamma < \gamma(\lambda)\}$ when we are in step~1 of the data generating algorithm in
Assumption~1, and to the event $R(\lambda)$ before the ball is drawn out of the urn, would usually
be classified as an objective probability. As a consequence it would usually be the case that the
similarity:
\begin{equation}
\label{equ24}
S(R(\lambda), \{\Gamma < \gamma(\lambda)\})
\end{equation}
is regarded as being very high, and hence, the probability $\lambda$ of the event
$\{\Gamma < \gamma(\lambda)\}$ is considered as being very strong.

In the definition of a univariate subjective fiducial distribution, i.e.\ Definition~2, the
probability of the event $\{\Gamma < \gamma(\lambda)\}$ is effectively treated as being $\lambda$
after the data have been observed. To determine what strength ought to be assigned to this
probability, let us consider a modified version of one of the abstract scenarios that were outlined
in Bowater~(2017b).

\pagebreak
In particular, suppose that someone, who will be referred to as the selector, randomly draws a ball
out of the urn that is associated with the set $R$ and then, without looking at the ball, hands it
to an assistant. The assistant, by contrast, looks at the ball, but conceals it from the selector,
and then places it under a cup. The selector believes that the assistant smiled when he looked at
the ball.

Under these conditions, the selector is asked to assign a probability to the event of the number on
the ball being less than or equal to $\lambda m$, where generally $\lambda$ can be any given value
in $\Lambda$, but it may be helpful to assume that $\lambda$ is not too close to 0 or 1. Let this
event be denoted as $R^*(\lambda)$. Finally, we will assume that it is uncertain whether the
assistant knew from the outset that the selector would be asked to assign a probability to this
particular event.

Clearly in this scenario, a smile by the assistant would in general need to be taken into account,
since it could imply that it is less likely or more likely that the event $R^*(\lambda)$ has taken
place. However, the selector may feel that, if the assistant had indeed smiled, he would not have
understood the smile's meaning. For this reason, he may decide that the probability that the number
on the ball is less than or equal to $\lambda m$ is what it was before the ball was drawn out of
the urn, i.e.\ it is $\lambda$.

Now, it would seem undeniable that this probability is a subjective probability as it depends on a
subjective judgement regarding the meaning of a supposed smile. However, given his lack of
understanding about this meaning, the selector may feel that the similarity
$S(R^*(\lambda),R(\lambda))$ is very high, and hence that the probability $\lambda$ of the event
$R^*(\lambda)$ is very strong or, to put it another way, that this probability can be regarded in
a certain sense as being almost like a physical probability, i.e.\ a probability derived simply
through considerations of physical symmetry.

In addressing the main issue of what strength should be assigned after the data have been observed
to the probability $\lambda$ of the event $\{\Gamma < \gamma(\lambda)\}$, an analogy can be drawn
between the supposed smile of the assistant in this abstract scenario and the event of observing
the data $x$. In particular, under the assumptions of Section~\ref{sec8}, if little or nothing was
known about the parameter $\theta_j$ before the data were observed, the event of observing the data
should have little or no meaning in terms of how it should affect the probability of the event
$\{\Gamma < \gamma(\lambda)\}$. Therefore, although the similarity in equation~(\ref{equ24}) when
judged after the data have been observed may be considered as being marginally less than what it
was before the data were observed, it still ought to be considered as being very high, and hence,
the probability $\lambda$ of the event $\{\Gamma < \gamma(\lambda)\}$ still ought to be regarded as
being very strong.

Obviously, if there had been strong beliefs about $\theta_j$ before the data were observed, then it
may be quite clear how observing the data $x$ should affect the probability of the event
$\{\Gamma < \gamma(\lambda)\}$.
For example, if such beliefs had been strong enough so that they could have been adequately
summarised by placing a probability distribution over $\theta_j$, then the probability of the event
in question $\{\Gamma < \gamma(\lambda)\}$ after the data have been observed could be determined by
using Bayes' theorem.
Nevertheless, the fiducial argument is traditionally applied under the assumption that little or
nothing was known about the parameters of interest before the data were observed, and the present
work will not deviate from this tradition.

Under this assumption and taking into account that the full conditional fiducial density of
$\theta_j$, i.e.\ the density $f(\theta_j\,|\,\theta_{-j}, x)$, is fully defined by the density
function of the primary r.v.\ $\Gamma$ and known constants, the probability $\lambda$ of the event
$\{\theta_j < \theta_j(\lambda)\}$ after the data have been observed, where $\theta_j(\lambda)$ is
defined by:
\vspace{2ex}
\[
\int_{-\infty}^{\theta_j(\lambda)} f(\theta_j\,|\,\theta_{-j},x) d\theta_j = \lambda
\vspace{2ex}
\]
and $\lambda$ is any given member of $\Lambda$, ought to be regarded as being a very strong
probability, or in other words, in spite of this probability being a fiducial \pagebreak
probability it ought to be regarded in a certain sense as being almost like a physical probability.

Let us now consider how strong we ought to regard the type of probabilities that are derived by
integrating over the joint fiducial density of all the parameters $\theta$ that is determined
according to the framework outlined in Sections~\ref{sec3} to~\ref{sec14}.
First, let us assume that the set of full conditional fiducial densities of the parameters $\theta$
that are determined using the method outlined in Section~\ref{sec8} are compatible and the joint
fiducial density of these parameters that they define is unique.
In this case, if probabilities derived by integrating over the full conditional fiducial density
$f(\theta_j\,|\,\theta_{-j},x)$ with respect to the parameter $\theta_j$ are considered as being
very strong for all $j \in \{1,2,\ldots,k\}$, then since the joint fiducial density of the
parameters $\theta$ is fully and directly defined by the set of full conditional densities in
question, it can be argued that the probabilities derived by integrating over this joint density
ought to be treated as though they are almost equivalent to physical probabilities.

On the other hand, if the full conditional fiducial densities of the parameters $\theta$ are
incompatible, then the same type of argument is still reasonably sound under the condition, which
is fair to expect would often be satisfied, that these densities are, nevertheless, close
approximations to the full conditional densities of the joint fiducial density of these parameters
that is determined in the most appropriate way by using a Gibbs sampler within the framework
outlined in Sections~\ref{sec3} to~\ref{sec14}.
Also, probabilities derived by integrating over any given one of these latter full conditional
densities could be regarded as actually being stronger than probabilities derived by integrating
over the full conditional fiducial density of the same parameter that was directly determined by
using the method described in Section~\ref{sec8}, since it is known that the former density
function belongs to a set of full conditional densities that have been adjusted so that they are
compatible.

\vspace{3ex}
\section{Comparing subjective fiducial inference to Bayesian inference}

As mentioned in the Introduction, given that in many cases the standard fiducial distribution is
equal to the posterior distribution obtained through Bayes' theorem for a given choice of the prior
distribution, it has become a convention to claim that, to a large extent, fiducial inference is
indistinguishable from Bayesian inference.
For this reason, it is worth comparing subjective fiducial inference to Bayesian inference. This
comparison will be carried out using the definition of probability upon which the present work is
based, i.e.\ the definition outlined in Bowater~(2017a) under the name of type B probability and in
Bowater~(2017b).
It is therefore necessary to apply this definition to the probabilities used in Bayes' theorem.

The probabilities that enter into Bayes' theorem are provided by the sampling density
$g(x\,|\,\theta)$ and the prior density $p(\theta)$.
Probabilities obtained by integrating over the sampling density can usually be regarded as physical
probabilities, and therefore, under the important assumption that the sampling model actually
generated the observed data, these probabilities can usually be considered as being very strong.
If probabilities obtained by integrating over the prior density $p(\theta)$ are also very strong,
it would seem logical to regard the probabilities produced by Bayes' theorem, i.e.\ those obtained
by integrating over the resulting posterior density $p(\theta\,|\,x)$, as also being very strong.
This would be the case, for example, if the prior distribution represents the uncertainty
concerning the outcome of a well-understood physical experiment.
Similarly, if the prior density of the parameter of interest $\theta$ is elicited on the basis of
the subjective opinion of one or a number of scientists or scholars who have quite detailed
knowledge about the relative plausibility of different values for this parameter, then
probabilities obtained by integrating over the prior density in question may well be considered as
being reasonably strong and, as a result, it is arguable that the posterior probabilities that
correspond to the use of this prior density should also be regarded as being reasonably strong.

However, to make a direct comparison with subjective fiducial inference, it needs to be assumed
that little or nothing was known about the model parameters $\theta$ before the data were observed.
If $\theta$ is a single parameter unrestricted on the real line, it is common to try to represent
this lack of knowledge by placing a diffuse symmetric prior density over $\theta$ centred at some
given value for its median, which will be denoted as the value $\theta^{*}(0.5)$.
Assuming that this has been done, let us consider the similarity between the event $R(\lambda)$ as
specified in the previous section and the event of $\theta$ being less than the value
$\theta^{*}(\lambda)$, i.e.\ the similarity $S(R(\lambda), \{\theta < \theta^{*}(\lambda)\})$,
where $\theta^{*}(\lambda)$ is defined by the expression:
\vspace{1.5ex}
\[
\int_{-\infty}^{\theta^{*}(\lambda)} p(\theta) d\theta = \lambda
\vspace{1.5ex}
\]
in which $p(\theta)$ is the chosen prior density and $\lambda \in \Lambda$.
Notice that no matter how diffuse the prior density is chosen to be, this similarity is likely to
be regarded as being very low for any given value of $\lambda$.
For example, if $\lambda=0.5$ then $S(R(\lambda), \{\theta < \theta^{*}(\lambda)\})$ is effectively
the similarity between the event of drawing a red ball out of an urn that contains an equal number
of red balls and blue balls, and the event of $\theta$ being less than $\theta^{*}(0.5)$, which is
clearly going to be a very low similarity as the choice of the median $\theta^{*}(0.5)$ is
completely arbitrary.
This implies that for any value $\lambda \in \Lambda$, the probability $\lambda$ of the event
$\{\theta < \theta^{*}(\lambda)\}$ is likely to be considered as being very weak.

Since posterior probabilities are derived through Bayes' theorem simply by reweighting prior
probabilities, that is by normalising the probabilities that are the result of multiplying prior
probabilities by the likelihood function, it would seem difficult to make the argument, on the
basis of the reasoning that underlies this theorem, that such probabilities should be regarded as
being generally stronger than prior probabilities. Of course, this does not discount the
possibility that this particular argument could be sometimes justified by using a non-Bayesian form
of reasoning, but this would raise the question as to whether Bayesian theory is actually being
applied.

\pagebreak
Therefore, if a proper prior density is used to try to represent a lack of knowledge about an
unrestricted parameter $\theta$ before the data are observed, then it could be argued that it would
be difficult to use Bayesian reasoning to make any kind of claim that probabilities not equalling
zero or one that are obtained by integrating over the resulting posterior density ought to be
regarded as being anything other than very weak probabilities.
A similar argument could be presented in the case where the parameter $\theta$ is restricted on the
real line.

This is a clear inadequacy of Bayesian inference. Also, notice that the severe criticisms that were
highlighted in Section~\ref{sec9} of what was referred to, in that section, as objective Bayesian
inference do not apply to the theory of subjective fiducial inference outlined in the present
paper.
Therefore, a strong case has been made that subjective fiducial inference is superior to Bayesian
inference in scenarios where little or nothing was known about the parameters of the model of
interest before the data were observed.

\vspace{3ex}
\section{Open issues}

In this closing section, we will briefly discuss some open issues concerning subjective fiducial
inference.

As pointed out in Lindley~(1958), standard fiducial inference can be incoherent in the sense that
treating a fiducial distribution for a given parameter that is derived on the basis of a data set
$x^{(1)}$ as a prior distribution in a Bayesian analysis of another independent data set $x^{(2)}$
does not lead to a posterior distribution that is equal to the fiducial distribution of the
parameter that corresponds to the combined data set $\{x^{(1)},x^{(2)}\}$.
Much has been made of the existence of this anomaly, however little attention has been given to its
practical consequences. In particular, little research has been done into establishing in what
situations fiducial inference fails to at least approximately satisfy this coherency condition,
especially when either or both of the samples $x^{(1)}$ and $x^{(2)}$ are at least moderately
sized.

Furthermore, in cases where there may be a substantial difference between the fiducial distribution
derived on the basis of the data set $\{x^{(1)},x^{(2)}\}$ and the posterior distribution that
results from analysing this data set by combining Bayesian and fiducial inference in the
aforementioned way, a sensible strategy exists for choosing between these two distributions.
In particular, if the first data set $x^{(1)}$ is large enough so that the fiducial distribution
that corresponds to this first data set is considered to be a very strong distribution, according
to the criteria given in Section~\ref{sec7} and in Bowater~(2017b), then it would seem sensible to
regard the posterior distribution that results from analysing the second data set $x^{(2)}$ with
this fiducial distribution as the prior distribution as providing the most appropriate inferences
about the parameters of interest on the basis of the data set $\{x^{(1)}, x^{(2)}\}$.
On the other hand, if the first data set $x^{(1)}$ is so small that some doubts exist with regard
to classifying the fiducial distribution that corresponds to this first data set as being a very
strong distribution then, in general, the most appropriate inferences about the parameters
concerned could be considered as being provided by the fiducial distribution that corresponds to
the combined data set $\{x^{(1)},x^{(2)}\}$.

As was illustrated in Section~\ref{sec10}, when applying subjective fiducial inference as defined
in the present paper to the univariate case as specified in Section~\ref{sec8}, it is not necessary
that the fiducial statistic is a sufficient statistic for the unknown parameter $\theta_j$, since
if a univariate sufficient statistic for $\theta_j$ does not exist, it can be defined to be any
one-to-one function of a unique maximum likelihood estimator of $\theta_j$.
This represents a departure from standard fiducial inference which, as is easy to see, will have
some consequences in terms of the coherency issue that has just been raised, but which nevertheless
substantially opens up the range of applications of the methodology that has been discussed.
It is left as an open issue as to whether and by how much subjective fiducial inference could
perform better in some cases if the fiducial statistic was allowed to be another type of
non-sufficient statistic.

Finally we note that, in general, it would appear that subjective fiducial inference is more
computationally demanding than Bayesian inference. However, on the basis of what has been seen in
recent years in relation to Bayesian inference, it is reasonable to anticipate that major advances
could be achieved with respect to the computational aspects of subjective fiducial inference, which
will gradually extend its range of applications, as well as improving the accuracy by which the
kind of functions on which it relies can be approximated.

\vspace{5ex}
\noindent
{\bf References}

\vspace{0.5ex}
\begin{description}

\setlength{\itemsep}{1.5ex}

\item[] Alary, D., Gollier, C. and Treich, N. (2013).\ The effect of ambiguity aversion on
insurance and self-protection.\ \emph{The Economic Journal}, {\bf 123}, 1188--1202.

\item[] Arnold, B. C., Castillo, E., Sarabia, J. M. (2002).\ Exact and near compatibility of
discrete conditional distributions.\ \emph{Computational Statistics and Data Analysis}, {\bf 40},
231--252.

\item[] Arnold, B. C. and Press, S. J. (1989).\ Compatible conditional distributions.\
\emph{Journal of the American Statistical Association}, {\bf 84}, 152--156.

\item[] Berger, J. O. and Bernardo, J. M. (1992).\ On the development of the reference prior
method.\ In \emph{Bayesian statistics 4}, Eds.\ J. M. Bernardo, J. O. Berger, A. P. Dawid and A. F.
M. Smith, pp.\ 35--49, Oxford University Press, London.

\item[] Bernardo, J. M. (1979).\ Reference posterior distributions for Bayesian inference (with
discussion).\ \emph{Journal of the Royal Statistical Society, Series B}, {\bf 41}, 113--147.

\item[] Bowater, R. J. (2017a).\ A formulation of the concept of probability based on the use of
experimental devices.\ \emph{Communications in Statistics:\ Theory and Methods}, {\bf 46},
4774--4790.

\setlength{\itemsep}{1ex}

\item[] Bowater, R. J. (2017b).\ A defence of subjective fiducial inference.\ \emph{AStA Advances
in Statistical Analysis}, {\bf 101}, 177--197.

\item[] Brooks, S. P. and Gelman, A. (1998).\ General methods for monitoring convergence of
iterative simulations.\ \emph{Journal of Computational and Graphical Statistics}, {\bf 7},
434--455.

\item[] Buehler, R. J. and Feddersen, A. P. (1963).\ Note on a conditional property of
Student's~t.\ \emph{Annals of Mathematical Statistics}, {\bf 34}, 1098--1100.

\item[] Chen, S-H. and Ip, E. H. (2015).\ Behaviour of the Gibbs sampler when conditional
distributions are potentially incompatible.\ \emph{Journal of Statistical Computation and
Simulation}, {\bf 85}, 3266--3275.

\item[] Chen, S-H., Ip, E. H. and Wang, Y. J. (2011).\ Gibbs ensembles for nearly compatible and
incompatible conditional models.\ \emph{Computational Statistics and Data Analysis}, {\bf 55},
1760--1769.

\item[] Cowles, M. K. and Carlin, B. P. (1996).\ Markov chain Monte Carlo convergence diagnostics:\
a comparative review.\ \emph{Journal of the American Statistical Association}, {\bf 91}, 883--904.

\item[] Dempster, A. P. (1968).\ A generalization of Bayesian inference (with discussion).\
\emph{Journal of the Royal Statistical Society, Series B}, {\bf 30}, 205--247.

\item[] Edwards, W., Lindman, H. and Savage, L. J. (1963).\ Bayesian statistical inference for
psychological research.\ \emph{Psychological Review}, {\bf 70}, 193--242.

\item[] Ellsberg, D. (1961).\ Risk, ambiguity and the Savage axioms.\ \emph{Quarterly Journal of
Economics}, {\bf 75}, 643--669.

\item[] de Finetti, B. (1974).\ \emph{Theory of Probability}, vol. 1., Wiley, Chichester.

\item[] de Finetti, B. (1975).\ \emph{Theory of Probability}, vol. 2., Wiley, Chichester.

\item[] Fisher, R. A. (1930).\ Inverse probability.\ \emph{Mathematical Proceedings of the
Cambridge Philosophical Society}, {\bf 26}, 528--535.

\item[] Fisher, R. A. (1935).\ The fiducial argument in statistical inference.\ \emph{Annals of
Eugenics}, {\bf 6}, 391--398.

\item[] Fisher, R. A. (1956).\ \emph{Statistical Methods and Scientific Inference}, 1st ed., Hafner
Press, New York [2nd ed., 1959; 3rd ed., 1973].

\item[] Fraser, D. A. S. (1966).\ Structural probability and a generalization.\ \emph{Biometrika},
{\bf 53}, 1--9.

\item[] Fraser, D. A. S. (1972).\ Events, information processing, and the structured model.\ In
\emph{Foundations of Statistical Inference}, Eds.\ V. P. Godambe and D. A. Sprott, pp. 32--55,
Holt, Renehart and Winston, Toronto.

\item[] Gelfand, A. E. and Smith, A. F. M. (1990).\ Sampling-based approaches to calculating
marginal densities.\ \emph{Journal of the American Statistical Association}, {\bf 85}, 398--409.

\item[] Gelman, A. and Rubin, D. B. (1992).\ Inference from iterative simulation using multiple
sequences.\ \emph{Statistical Science}, {\bf 7}, 457--472.

\item[] Geman, S. and Geman, D. (1984).\ Stochastic relaxation, Gibbs distributions and the
Bayesian restoration of images.\ \emph{IEEE Transactions on Pattern Analysis and Machine
Intelligence}, {\bf 6}, 721--741.

\item[] Ghosh, M. (2011).\ Objective priors:\ an introduction for frequentists (with discussion).\
\emph{Statistical Science}, {\bf 26}, 187--211.

\item[] Gilboa, I. and Schmeidler, D. (1989).\ Maximin expected utility with non-unique prior.\
\emph{Journal of Mathematical Economics}, {\bf 18}, 141--153.

\item[] Goutis, C. and Casella, G. (1995).\ Frequentist post-data inference.\ \emph{International
Statistical Review}, {\bf 63}, 325--344.

\item[] Hannig, J. (2009).\ On generalized fiducial inference.\ \emph{Statistica Sinica}, {\bf 19},
491--544.

\item[] Hannig, J., Iyer, H., Lai, R. C. S. and Lee, T. C. M. (2016).\ Generalized fiducial
inference:\ a review and new results.\ \emph{Journal of the American Statistical Association},
{\bf 111}, 1346--1361.

\item[] Jeffreys, H. (1961).\ \emph{Theory of Probability}, 3rd ed., Oxford University Press,
Oxford.

\item[] Kass, R. E. and Wasserman, L. (1996).\ The selection of prior distributions by formal
rules.\ \emph{Journal of the American Statistical Association}, {\bf 91}, 1343--1370.

\item[] Kuo, K-L., Song, C-C. and Jiang, T. J. (2017).\ Exactly and almost compatible joint
distributions for high-dimensional discrete conditional distributions.\ \emph{Journal of
Multivariate Analysis}, {\bf 157}, 115--123.

\item[] Kuo, K-L. and Wang, Y. J. (2011).\ A simple algorithm for checking compatibility among
discrete conditional distributions.\ \emph{Computational Statistics and Data Analysis}, {\bf 55},
2457--2462.

\item[] Lecoutre, B. and Poitevineau, J. (2014).\ \emph{The Significance Test Controversy
Revisited:\ the Fiducial Bayesian Alternative}, Springer, Heidelberg.

\item[] Lindley, D. V. (1958).\ Fiducial distributions and Bayes' theorem.\ \emph{Journal of the
Royal Statistical Society, Series B}, {\bf 20}, 102--107.

\item[] Lindley, D. (1997).\ Some comments on `Non-informative priors do not exist'.\ \emph{Journal
of Statistical Planning and Inference}, {\bf 65}, 182--184.

\item[] Liu, C. and Martin, R. (2015).\ Frameworks for prior-free posterior probabilistic
inference.\ \emph{WIREs Computational Statistics}, {\bf 7}, 77--85.

\item[] Martin, R. and Liu, C. (2015).\ \emph{Inferential Models:\ Reasoning with Uncertainty},
CRC Press, Boca Raton.

\item[] Mur\'e, J. (2019).\ Optimal compromise between incompatible conditional probability
distributions, with application to objective Bayesian kriging.\ \emph{ESAIM:\ Probability and
Statistics}, {\bf 23}, 271--309.

\item[] Savage, L. J. (1954).\ \emph{The Foundations of Statistics}, Wiley, New York.

\item[] Seidenfeld, T. (1979).\ Why I am not an objective Bayesian; some reflections prompted by
Rosenkrantz.\ \emph{Theory and Decision}, {\bf 11}, 413--440.

\item[] Shafer, G. (1976).\ \emph{A Mathematical Theory of Evidence}, Princeton University Press,
Princeton.

\item[] Spetzler, C. S. and Stael von Holstein, C. A. S. (1975).\ Probability encoding in decision
analysis.\ \emph{Management Science}, {\bf 22}, 340--358.

\item[] Wilkinson, G. N. (1977).\ On resolving the controversy in statistical inference (with
discussion).\ \emph{Journal of the Royal Statistical Society, Series B}, {\bf 39}, 119--171.

\item[] Xie, M. and Singh, K. (2013).\ Confidence distribution, the frequentist distribution
estimator of a parameter:\ a review.\ \emph{International Statistical Review}, {\bf 81}, 3--39.

\end{description}

\end{document}